\documentclass[11pt]{article}
\usepackage{CJK}
\usepackage{ws-rotating}
\usepackage{epstopdf}
\usepackage{epsfig}
\usepackage{color}
\usepackage{amsmath}
\usepackage{amssymb}
\usepackage{paralist}
\usepackage{graphicx}
\usepackage{epstopdf}
\usepackage[numbers,sort&compress]{natbib}
\usepackage{float}

\newtheorem{mythm}{\bf{Theorem}\rm}
\newtheorem{myprop}{\bf{Proposition}\rm}
\newtheorem{mylem}{\bf{Lemma}\rm}

\newtheorem{remark}{Remark}

\begin{document}
\pagenumbering{arabic}

\title{Does solitary wave solution persist for the long wave equation with small perturbations?
}
\author{Hang Zheng$^{1,2,3}$\,\,\,\, Y-H. Xia$^{1}$\footnote{Corresponding author: Y-H. Xia,
 Email: xiadoc@outlook.com; yhxia@zjnu.cn.}
  \\
{\small \em 1. Department of Mathematics, Zhejiang Normal University, Jinhua, 321004, China}\\
{\small   \,\,\,  yhxia@zjnu.cn; xiadoc@outlook.com}\\
{\small \em 2.  Department of Mathematics and Computer, Wuyi University, Wuyishan, 354300, China}\\
{\small \em 3.  Fujian key Laboratory of Big Data Application and Intellectualization for Tea Industry}\\
{\small  \,\,\,  zhenghang513@zjnu.edu.cn; zhenghwyxy@163.com}\\
}

\date{}
\maketitle \pagestyle{myheadings} \markboth{}{}

\noindent

\begin{abstract}
{ Persistence of solitary wave solutions of the regularized long wave equation with small perturbations
are investigated by the geometric singular perturbation theory. Two different kinds of the perturbations are considered in this paper: one is the weak backward diffusion and dissipation, the other is the Marangoni effects.
%To prove the persistence of solitary wave solution, the analytical expression of Melnikov integral for the perturbed RLW equation is directly obtained. %which is different from the previous literature.
 Indeed, the solitary wave persists under small perturbations. Furthermore, the different perturbations do affect the  proper wave speed $c$ ensuring the persistence of the solitary waves. Finally, numerical simulations are utilized to confirm the theoretical results.
  %It is proved that the solitary wave solution persists with an appropriate speed $c$. Interestingly, comparing with BDD perturbation, the term $(uu_x)_x$ in Marangoni effects perturbation induces a different speed $\widetilde{c}$ for the existence of solitary wave solution.
  }
\end{abstract}


\section{Introduction}

\subsection{History of the perturbed long wave equation}
\subsubsection{Unperturbed long wave equation }
 { Shallow water models have been one of the research hotspots of nonlinear partial differential equations.
  %because it models many nonlinear complex phenomena in fields such as physics, mechanics and biology, etc.
  There are many shallow water models appearing in different forms, such as KdV equation (Korteweg and de Vries \cite{KdV-PM}), CH equation (Camassa and Holm \cite{CH-PRE}%, Constantin\cite{Constantin-JFA,Constantin-RSPA}, Constantin and Escher \cite{Constantin-acta}
  ), DP equation (Degasperis and Procesi \cite{DP-NJ}), BBM equation (Benjamin et al.\cite{BBM-PTRSLA}), GN equation (Green and Naghdi \cite{GN-JFM}), %shallow water models with Coriolis effect (see e.g., Gui et al. \cite{Gui1,Gui2} and Chen et al. \cite{Gui3}), some other forms of shallow water models (see e.g., Chu et al. \cite{CHUjf-JDE2020,CHUjf-SPAM}, Roberts et al. \cite{RAJ-SIAM}, Fu et al. \cite{fy-JFA}), 
  and so on. The classical KdV equation is formulated as follows
\begin{equation}\label{KdV}
u_t+uu_x+u_{xxx}=0.
\end{equation}
%After that, a lot of generalized KdV-type equations have been proposed.
}In 1967, Peregrine \cite{PHD-JFM1,PHD-JFM2} put forward a regularized long-wave (RLW, for short) equation which regarded as a alternative model to KdV equation of the form
\begin{equation}\label{RLW}
u_t+u_x+uu_x-u_{xxt}=0.
\end{equation}
The equation is of many attractive features which can be concluded as follows (see \cite{BMR-JNS}): (i) more closely matching that of the full Euler equations to describe the two-dimensional motion of the free-surface oscillations of an ideal liquid under the force of gravity, (ii) it is better to integrate numerically. Later, to describe surface water waves in a uniform channel, a regularized form of the KdV equation was studied by Benjamin et al. (see \cite{BBM-PTRSLA}), and the BBM equation is given by
\begin{equation}\label{BBM}
u_t+uu_x-u_{xxt}=0.
\end{equation}
\subsubsection{KS perturbation}
However, in a more realistic situation, it is unavoidable that the existence of uncertainty or perturbations. Thus, small perturbation terms should be considered in modelling the problems. {For shallow water wave equation, weak backward diffusion $u_{xx}$ and dissipation $u_{xxxx}$ are called Kuramoto-Sivashinsky (KS, for short) which are one of common perturbations.} To explain the wave motions on
a liquid layer over an inclined plane, Topper and Kawahara \cite{TTK-Jpsj} analyzed a PDE as follows:
\begin{equation}\label{KBK}
u_t+uu_x+\alpha u_{xx}+\beta u_{xxx}+\gamma u_{xxxx}=0,
\end{equation}
where $\alpha$, $\beta$ and $\gamma$ are positive parameters. If the inclined plane and surface tension are relatively long and weak, it implies $u_{xx}$ and $u_{xxxx}$ are real small terms such that they are considered to be some small perturbations. Then Eq. (\ref{KBK}) is transformed to the following equation
\begin{equation}\label{KBK-p}
u_t+uu_x+ u_{xxx}+\varepsilon( u_{xx}+u_{xxxx})=0,
\end{equation}
where $0<\varepsilon\ll1$, Eq. (\ref{KBK-p}) is usually called as a perturbed KdV equation. %For convenience, we abbreviate backward diffusion $u_{xx}$ and dissipation $u_{xxxx}$ to BDD perturbation.
 There are many monographs studying the existence of the traveling wave solution under KS perturbation. For instance, Derks and Gils \cite{DG-JJIAM}, Ogawa \cite{OT-HMJ} investigated the existence of traveling wave solution of Eq. (\ref{KBK-p}). Yan et al. \cite{YW-MMA} considered a generalized KdV equation with KS perturbation and proved the existence of solitary wave and periodic wave solution. Chen et al. \cite{CAY-AML,CAY-JDE,CAY-QTDS} explained the existence of traveling wave solution of a generalized KdV equation and BBM equation with KS perturbation.
\subsubsection{Marangoni effect perturbation}
In addition to KS perturbation, Marangoni effect (ME) is { an alternative perturbation in modelling,} see \cite{gpl-PF}. Generally, the change of surface tension gradient at the interface of two phases {leads to} the occurrence of ME. In other words, ME presents an additional nonlinear term $(uu_x)_x$ when it is considered on the surface of a thin layer. Velarde et al. \cite{VMG-NY} proposed the following equation
\begin{equation}\label{ME}
u_t+2\alpha_1uu_x+\alpha_2 u_{xx}+\alpha_3 u_{xxx}+\alpha_4 u_{xxxx}+\alpha_5(uu_x)_x=0
\end{equation}
which involving ME into the one-way long-wave assumption. And nonlinear term $(uu_x)_x$ generates from ME, which describes the opposite to the B\'{e}nard convection. Sun and Yu \cite{SY2019} proved the existence of periodic wave and solitary wave solution of BBM equation under KS and generalized ME perturbations, the equation is described by
\begin{equation}\label{BBM-ME}
(u^3)_t+(u^4)_x+u_{xxx}+\epsilon E_i=0, \ i=1,2,
\end{equation}
where $E_1=u_{xx}+u_{xxxx}$ and $E_2=((\alpha_0+\alpha_1 u+\alpha_2 u^2)u_x)_x$. It shows that traveling wave are very sensitive to KS and weak ME perturbations.
{
\subsubsection{Model formulation}
}

In recent years, more and more papers concern the influence of distributed delay in the evolution of equations. Delays {are} more closely to show realistic phenomenon than the usual kind of models without delay. In 1989, Britton \cite{BNF-JTB} derived a model for a single biological population of the form
\begin{equation}\label{BP}
u_t=u(1+\alpha u-(1+\alpha)(f\ast u))+\bigtriangleup u,
\end{equation}
where $f\ast u$ is a local convolution in the spatial-temporal and defined by
\begin{equation}\label{defination}
(f\ast u)(x,t)=\int^t_{-\infty}f(t-s)u(x,s)ds.
\end{equation}
The kernel $f(t)$ satisfies $f:[0,\infty)\rightarrow [0,\infty)$ and
$$\int^{\infty}_0f(t)dt=1, \ \ tf(t)\in L^1((0,\infty),\mathbb{R}).$$
Here, we suppose that the kernel $f(t)$ is the strong generic delay kernel given by
\begin{equation}\label{function}
f(t)=\frac{4t}{\tau^2}e^{-{\frac{2t}{\tau}}}
\end{equation}
with averaging delay $\tau=\displaystyle \int^{\infty}_0tf(t)dt$. In this paper, we study a delayed RLW equation with  two different kinds of perturbations (KS and ME) as follows
\begin{equation}\label{PRLW}
u_t+u_x+(f\ast u)u_x-u_{xxt}+\tau P_i=0,\ \ (i=1,2),
\end{equation}
where $0<\tau\ll 1$ (a sufficiently small parameter), $P_1=u_{xx}+u_{xxxx}$ and $P_2=u_{xx}+(uu_x)_x+u_{xxxx}$. $P_1$ is the KS perturbation and $P_2$ represents the ME perturbation.
$f\ast u$ is a local convolution in the spatial-temporal, which is defined in (\ref{defination}) and (\ref{function}). It should be noted that BBM equation can be derived from equation (\ref{PRLW}). Most of the articles concern the traveling wave solution of unperturbed RLW equation (see e.g., Bona et al. \cite{BMR-JNS}, Dag et al. \cite{Dag-Saka-AMC}, Gardner et al. \cite{Gardner-CNME2}). Sun and Yu \cite{SY2019} considered the generalized ME perturbation in BBM equation {by Chebyshev criteria and Abelian integral}. Note that it is actually a regular perturbation when perturbed equation transformed into its travelling wave system. However, there is no paper considered the perturbed RLW equation with two different perturbations and delay by singular perturbation theory (GSP, for short). This paper aims to study the persistence of solitary wave solution of \eqref{PRLW} by GSP theory.

{\subsection{Motivation and contributions}
{
%\subsubsection{Application of singular perturbation theory to the wave equations}

It is not difficult to show the existence of a solitary wave solution to an unperturbed RLW equation. Does the solitary wave solution persist under small perturbations?
In this paper, we pay particular attention to answer the question.}
When a perturbed equation is a singularly perturbed system,
the geometric singular perturbation  theory (Fenichel \cite{FN-JDE} and Jones \cite{JCK-berlin}) is an effective tool to prove the persistence of traveling wave solution. There exist a lot of articles on this issue including perturbed KdV equation (see e.g., Chen et al. \cite{CAY-AML},  Derks and Gils \cite{DG-JJIAM}, Ogama \cite{OT-HMJ}), perturbed CH equation (see e.g., Du et al. \cite{DLL2018,DZ-JDE}, Qiu et al. \cite{QHM-CNSNS}), random systems (see e.g., Li et al. \cite{LJ-TAMS,LJ-DCDS,LJ-JDE}), perturbed BBM equation (see e.g., Chen et al. \cite{CAY-JDE}, Sun and Yu \cite{SY2019}), Fitzhugh-Nagumo equation (see Liu et al. \cite{Liuweishi-JDE1}), biological model (see e.g., Chen and Zhang \cite{Zhangxiang-SpM}, Wang and Zhang \cite{Zhangxiang-JDE}), perturbed sine-Gordon equation (see Derks et al. \cite{Derks}) and so on.} Du and Qiao \cite{DQ2020} utilized GSP theory, Fredholm orthogonality and asymptotic theory to prove the existence of traveling wave solutions in a Belousov-Zhabotinskii system with delay. Combining Poincar\'{e}-Bendixson theorem, Du et al. \cite{DZ-JDE2} also confirmed the persistence of solitary wave solution of a generalized Keller-Segel system with small cell diffusion.

{ Usually, it is necessary to compute the Melnikov function by using the GSP theory. For the perturbed Hamiltonian system, Melnikov function is associated to the Abelian integral.
It is an effective skill to prove existence of traveling wave solutions by detecting monotonicity of the ratio of Abelian integral (MRAI).} For example, Derks and Gils \cite{DG-JJIAM}, Ogawa \cite{OT-HMJ}, Chen et al. \cite{CAY-JDE,CAY-AML,CAY-QTDS}, Sun et al. \cite{SunNA,SY2019} computed the MRAI of corresponding perturbed equations to prove the existence of traveling wave solution. Besides, another effective way is to directly compute the exact representation of Melnikov function. Qiu et al. \cite{QHM-CNSNS}, Zhu et al. \cite{ZK-ND}, Cheng and Li \cite{cheng-dcds} computed the explicit expression of Melnikov function to discuss the existence of traveling wave solution in a perturbed generalized BBM equation, generalized CH equation delayed DP equation, respectively.
%Combining bifurcation theory of limit cycles (\cite{CHB-JMAA,Hanmaoan-book, Hanmaoan-chinese3}), the existence of the solitary wave solutions of the perturbed mKdV and mK$(3,1)$ equations are demonstrated by Zhang et al. \cite{Zhanglijun-ND,Zhanglijun-IJBC}.

{
Motivated by aforementioned works, we use the singular perturbation theory to study the delayed RLW equation \eqref{PRLW} with two different kinds of perturbations (KS and ME). We prove that the solitary wave solution does persist under small perturbations.
To this end, we proceed with three steps.

\noindent {\bf Step 1.} The exact expression of solitary wave solution of unperturbed RLW equation are obtained by using dynamical system method.

\noindent {\bf Step 2.} The homoclinic orbits are constructing by tracking invariant manifolds of corresponding ODEs such that the singular perturbed system reduces to a regular perturbed system.

\noindent {\bf Step 3.} The analytical expression of Melnikov integral is given for perturbed RLW equation.	Some elaborate analytical skills are employed to discuss the simple zero of the Melnikov integral associated to the perturbed RLW equation, due to its complexity. Finally, by Melnikov method and bifurcation theorem, we conclude that the solitary wave solution persists at proper speed.

%The highlights of this paper are summarized as follows:

%(i) GSP theory plays an essential role in dealing with singular perturbations and establishing the persistence of solitary wave solution.

%(ii) The addition of convolution incorporating KS and ME perturbations have more physical significance for modelling nonlinear phenomena. It should be noted that BBM equation can be derived from equation (\ref{PRLW}).

%(iii) Different from aforementioned works (see, e.g., \cite{DQ2020,DZ-JDE2,OT-HMJ,CAY-JDE,CAY-AML,CAY-QTDS,SunNA,SY2019}). In particular, Sun and Yu \cite{SY2019} considered the generalized ME perturbation in BBM equation, it is actually a regular perturbation when perturbed equation transformed into ODEs. However, the ME perturbation in ours work is a singular perturbation in corresponding ODEs.

%(iv) Different from previous literature, we directly compute the analytical expression of Melnikov integral of perturbed RLW equation to prove the persistence of solitary wave solution. Numerical simulations and comparison are utilized to verify the propositions and theorems.

\subsection{Outline of the paper}
We present some preliminary concept and lemmas in next section. We obtain the exact expression of solitary wave solution of unperturbed RLW equation in section 3. In section 4, we construct the homoclinic orbits and compute the Melnikov integral associated to the perturbed RLW equations. Then after analysis of the Melnikov integrals, we present our main results in this section. Finally, we present the numerical simulations to verify  the theoretical results.
}

%\subsection{The outline of the paper}
%The outline of the paper is as follows. Some preliminaries including GSP theory are introduced in section 2. The bifurcations of phase portraits and solitary wave solution of the unperturbed RLW equation are obtained in section 3. Based on the GSP approach, Melnikov method and bifurcation theorem, the persistence of solitary wave solution of perturbed RLW with KS and ME perturbations are respectively discussed in section 4. Numerical simulations and comparison are applied to verify previous theoretical results in section 5.

\section{Preliminaries}
{In this section, we introduce some known results which will be used in our proof.}
\begin{mylem}\label{Lemma1.1} (see \cite{FN-JDE})
Consider the system
\begin{equation}\label{System1}
\begin{cases}
x_1'=f(x_1,x_2,\epsilon),\\
x_2'=\epsilon g(x_1,x_2,\epsilon),
\end{cases}
\end{equation}
where $'=\frac{d}{dt}$, $x_1\in \mathbb{R}^n$, $x_2\in \mathbb{R}^l$ and $\epsilon$ is a real parameter and positive, the function $f$ and $g$ are $C^\infty$ on a set $U\times I$ where $U\subseteq \mathbb{R}^{n+l}$ is open, and $I$ is an open interval, containing $0$. If $\epsilon=0$, the system has a compact, normally hyperbolic manifold of points $M_0$ which is contained in the set $\{f(x_1,x_2,\epsilon)=0\}$. Then for sufficiently small positive $\epsilon$ and any $0<r<+\infty$. There exists a manifold $M_\epsilon$ such that following hold:

(i) which is locally invariant under the flow of system (\ref{System1});

(ii) which is $C^r$ in $x_1$, $x_2$ and $\epsilon$;

(iii) $M_\epsilon=\{(x_1,x_2):x_1=h^\epsilon(x_2)\}$ for some $C^r$ function $h^\epsilon(x_2)$ and $x_2$ in some compact $K$;

(iv) there exist locally invariant stable and unstable manifold $W^s(M_\epsilon)$ and $W^u(M_\epsilon)$ lying within $O(\epsilon)$ and being $C^r$ diffeomorphic to the stable and unstable manifold $W^s(M_0)$ and $W^u(M_0)$.
\end{mylem}

\begin{mylem}\label{Lemma1.2}
If $\tau\rightarrow 0$, then $(f\ast u)(x,t)\rightarrow u(x,t)$.
\end{mylem}
\emph{Proof. } According to expressions (\ref{defination}) and (\ref{function}). The $(f\ast u)(x,t)$ is defined by
\begin{equation}\label{function1}
\begin{array}{rcl}
(f\ast u)(x,t)= \displaystyle\int^{t}_{-\infty} f(t-s)u(x,s)ds
=\displaystyle\int^{t}_{-\infty} \frac{4(t-s)}{\tau^2}e^{-{\frac{2(t-s)}{\tau}}}u(x,s)ds.
\end{array}
\end{equation}
Let $\lambda=-\frac{2(t-s)}{\tau}$, it yields
\begin{equation}\label{function2}
\begin{array}{rcl}
(f\ast u)(x,t)=-\displaystyle\int^{0}_{-\infty}\lambda e^\lambda u(x,t+\frac{\lambda\tau}{2})d\lambda.
\end{array}
\end{equation}
Hence, we have
\begin{equation}\label{function3}
\begin{array}{rcl}
\lim\limits_{\tau\rightarrow0}(f\ast u)(x,t)&=&-\lim\limits_{\tau\rightarrow0}\displaystyle\int^{0}_{-\infty}\lambda e^\lambda u(x,t+\frac{\lambda\tau}{2})d\lambda\\
&=&-\displaystyle\int^{0}_{-\infty}\lambda e^\lambda u(x,t)d\lambda\\
&=&-u(x,t)\displaystyle\int^{0}_{-\infty}\lambda e^\lambda d\lambda\\
&=&u(x,t).
\end{array}
\end{equation}
The proof is completed.

\section{Solitary wave solution of the unperturbed RLW equation}
In this section, the solitary wave solution of unperturbed equation (\ref{PRLW})$|_{\tau=0}$ is obtained by dynamical system method.

We first introduce the following transformations:
\begin{equation}\label{transformations1}
u(x,t)=\phi(\xi),\ \ \ \xi=x-ct,
\end{equation}
then, equation (\ref{PRLW})$|_{\tau=0}$ is converted to an ordinary differential equation given by
\begin{equation}\label{ODEs}
(1-c)\phi'+\phi\phi'+c\phi'''=0,
\end{equation}
where $'=\frac{d}{d\xi}$. Integrating both sides of equation (\ref{ODEs}) once with respect to $\xi$, we have
\begin{equation}\label{ODEs1}
(1-c)\phi+\frac{1}{2}\phi^2+c\phi''=g,
\end{equation}
where $g$ is the integration constant ($g\in \mathbb{R}$). According to actual physical meanings, we are only interested in the case $c>0$.

Next, the bifurcation and phase portraits is studied by dynamical system method.
{\subsection{Homoclinic orbit of the travelling wave system \eqref{ODEs1}} }
Eq. (\ref{ODEs1}) is equivalent to the following planar system:
\begin{equation}\label{unperturbed system}
\begin{cases}
\phi'=y,\\
y'=\displaystyle\frac{1}{c}\big[(c-1)\phi-\frac{1}{2}\phi^{2}+g\big].
\end{cases}
\end{equation}
It is easy to obtain the corresponding first integral
\begin{equation}\label{First integral}
\begin{array}{rcl}
H(\phi,y)=\frac{1}{2}y^2-\frac{1}{c}[\frac{(c-1)}{2}\phi^2-\frac{1}{6}\phi^{3}+g\phi]=h.
\end{array}
\end{equation}
Denote that
\begin{equation}\label{expansion1}
\begin{array}{rcl}
&&\Delta=(c-1)^2+2g,\  H(\phi_1,0)=h_1,\  H(\phi_2,0)=h_2,\\
&&\phi_1=c-1-\sqrt{\Delta},\  \phi_2=c-1+\sqrt{\Delta}. \\
\end{array}
\end{equation}
Obviously, system (\ref{unperturbed system}) has two equilibrium $(\phi_1,0)$ and $(\phi_2,0)$ when $\Delta>0$. Here, we suppose that $c>0$ in the sense of physics significance.

\begin{myprop}\label{proposition}
If $(c-1)^2>-2g$, system (\ref{unperturbed system}) has a homoclinic orbit $\Gamma$ surrounding the center point $(\phi_2,0)$  to one saddle point $(\phi_1,0)$ when level curves defined by $H(\phi,y)=h_1$. (see Fig. 1).
\end{myprop}
\emph{Proof. } The linearized matrix of system (\ref{unperturbed system}) is
\begin{equation}\label{martix}
A=\left(
\begin{array}{ccccc}
0 & 1 \\
\displaystyle\frac{c-1-\phi}{c} & 0  \\
\end{array}
\right).
\end{equation}
Thus, if $(c-1)^2>-2g$,
\begin{equation}\label{phi1}
J(\phi_1,0)=\textmd{det} A_{(\phi_1,0)}=\left|
\begin{array}{ccccc}
0 & 1 \\
\frac{\sqrt{\Delta}}{c} & 0  \\
\end{array}
\right|=-\frac{\sqrt{\Delta}}{c}<0,\ \ \textmd{Trace}(A)_{(\phi_1,0)}=0,
\end{equation}
and
\begin{equation}\label{phi2}
J(\phi_2,0)=\textmd{det} A_{(\phi_2,0)}=\left|
\begin{array}{ccccc}
0 & 1 \\
-\frac{\sqrt{\Delta}}{c} & 0  \\
\end{array}
\right|=\frac{\sqrt{\Delta}}{c}>0,\ \ \textmd{Trace}(A)_{(\phi_2,0)}=0.
\end{equation}
{By the bifurcation theory of planar dynamical systems (see \cite{Li-Book,CHB-JMAA})}, we know that $(\phi_1,0)$ is a saddle point and $(\phi_2,0)$ is a center point.
Then, if the level curves defined by $H(\phi,y)=h_1$, there exists a homoclinic orbit $\Gamma$ (blue curve) surrounding the center point $(\phi_2,0)$ to one saddle point $(\phi_1,0)$ (see Fig. 1).

This completes the proof.
\begin{center}
\begin{tabular}{cc}
\epsfxsize=6.5cm \epsfysize=6.5cm \epsffile{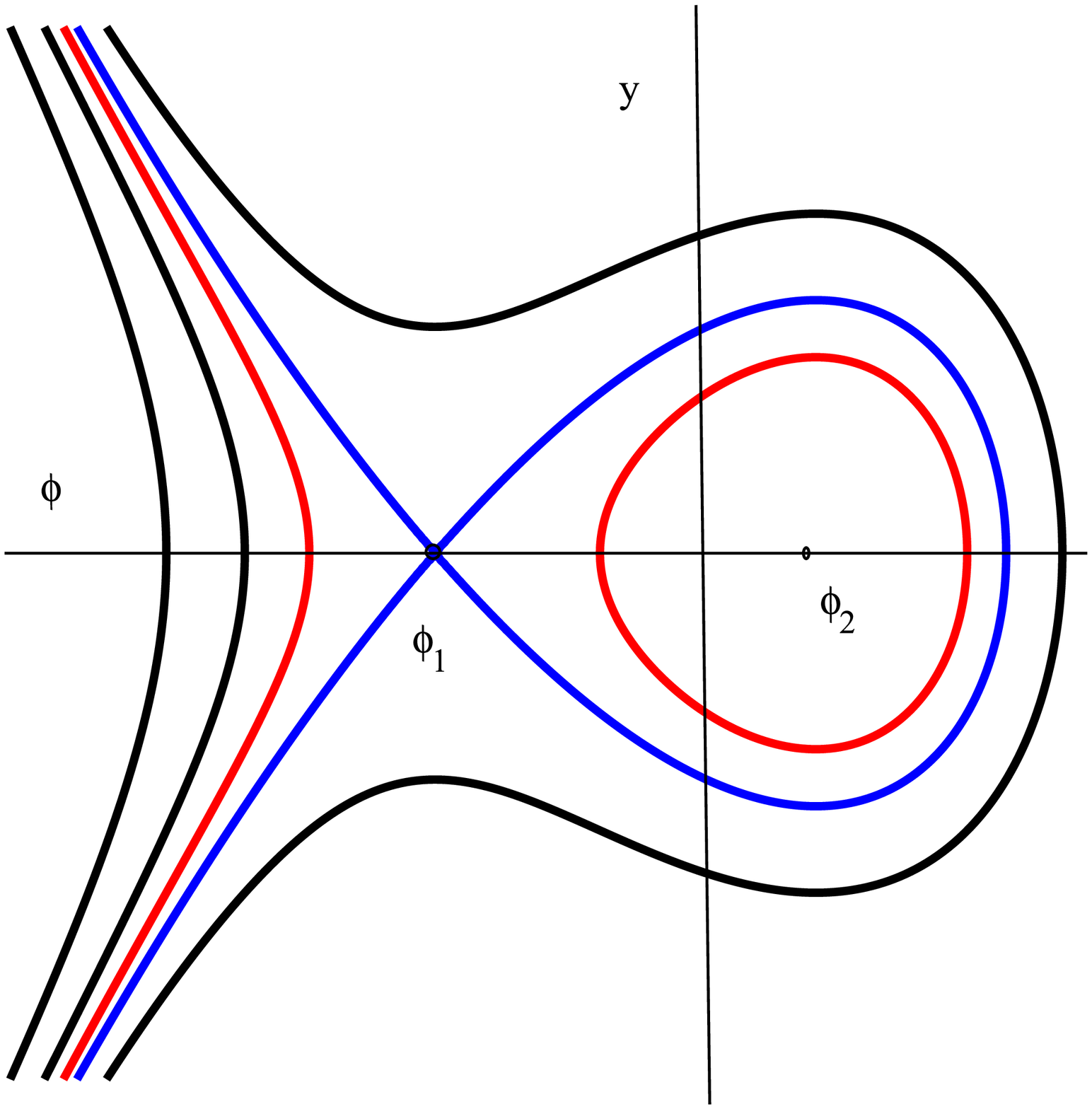}&
\epsfxsize=6.5cm \epsfysize=6.5cm \epsffile{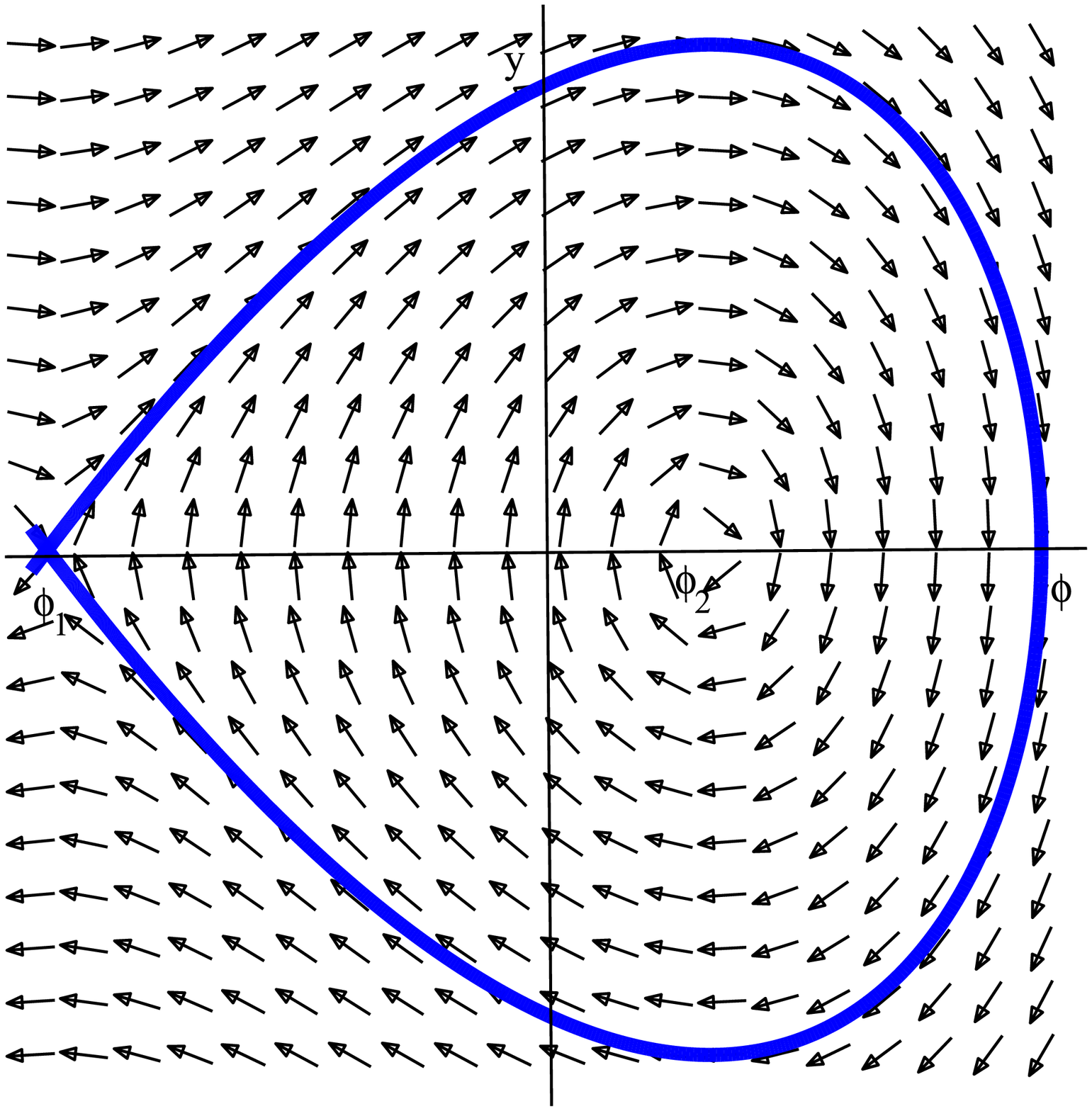}\\
\footnotesize{(a) Phase portraits  } & \footnotesize{ (b)  Homoclinic orbit }
\end{tabular}
\end{center}
\begin{center}
\footnotesize {{Fig. 1} \ \ Bifurcation of phase portraits of system (\ref{unperturbed system}) when $(c-1)^2>-2g$.}
\end{center}
{\subsection{Solitary wave solution of the unperturbed RLW equation} }
According to \textbf{Proposition} \ref{proposition} and (\ref{First integral}), the expression of $y$ can be written as
\begin{equation}\label{expression1}
\begin{array}{rcl}
y&=&\pm \sqrt{\frac{1}{3c}\phi^3+\frac{c-1}{c}\phi^2+\frac{2g}{c}\phi+\frac{2h_1}{c}}\\
 &=&\pm \sqrt{\frac{1}{3c}(\phi_r-\phi)(\phi-\phi_1)^2},
\end{array}
\end{equation}
where $\phi_r$ represents the right intersection of the homoclinic orbit with the $\phi$ axis and $\phi_r=2\sqrt{\Delta}+c-1$, it leads to
\begin{equation}\label{expression2}
\frac{d\phi}{d\xi}=\pm \sqrt{\frac{1}{3c}(\phi_r-\phi)(\phi-\phi_1)^2},
\end{equation}
which implies
\begin{equation}\label{expression3}
\xi=\int_{\phi_r}^{\phi}\frac{1}{\sqrt{\frac{1}{3c}(\phi_r-s)(s-\phi_1)^2}}ds.
\end{equation}
Then, the parametric representation of homoclinic orbit (see Fig. 2(a)) can be given by
\begin{equation}\label{homoclinic1}
\begin{array}{rcl}
\phi(\xi)=\phi_r-(\phi_r-\phi_1)\textmd{tanh}^2(\frac{1}{2}\sqrt{\frac{\phi_r-\phi_1}{3c}}\xi).
\end{array}
\end{equation}
Therefore, it corresponds to a bright solitary wave solution of Eq. (\ref{PRLW})$|_{\tau=0}$ (see Fig. 2(b)):
\begin{equation}\label{homoclinic2}
\begin{array}{rcl}
u(x,t)=\phi_r-(\phi_r-\phi_1)\textmd{tanh}^2(\frac{1}{2}\sqrt{\frac{\phi_r-\phi_1}{3c}}(x-ct)).
\end{array}
\end{equation}
\begin{center}
\begin{tabular}{cc}
\epsfxsize=6.5cm \epsfysize=6.5cm \epsffile{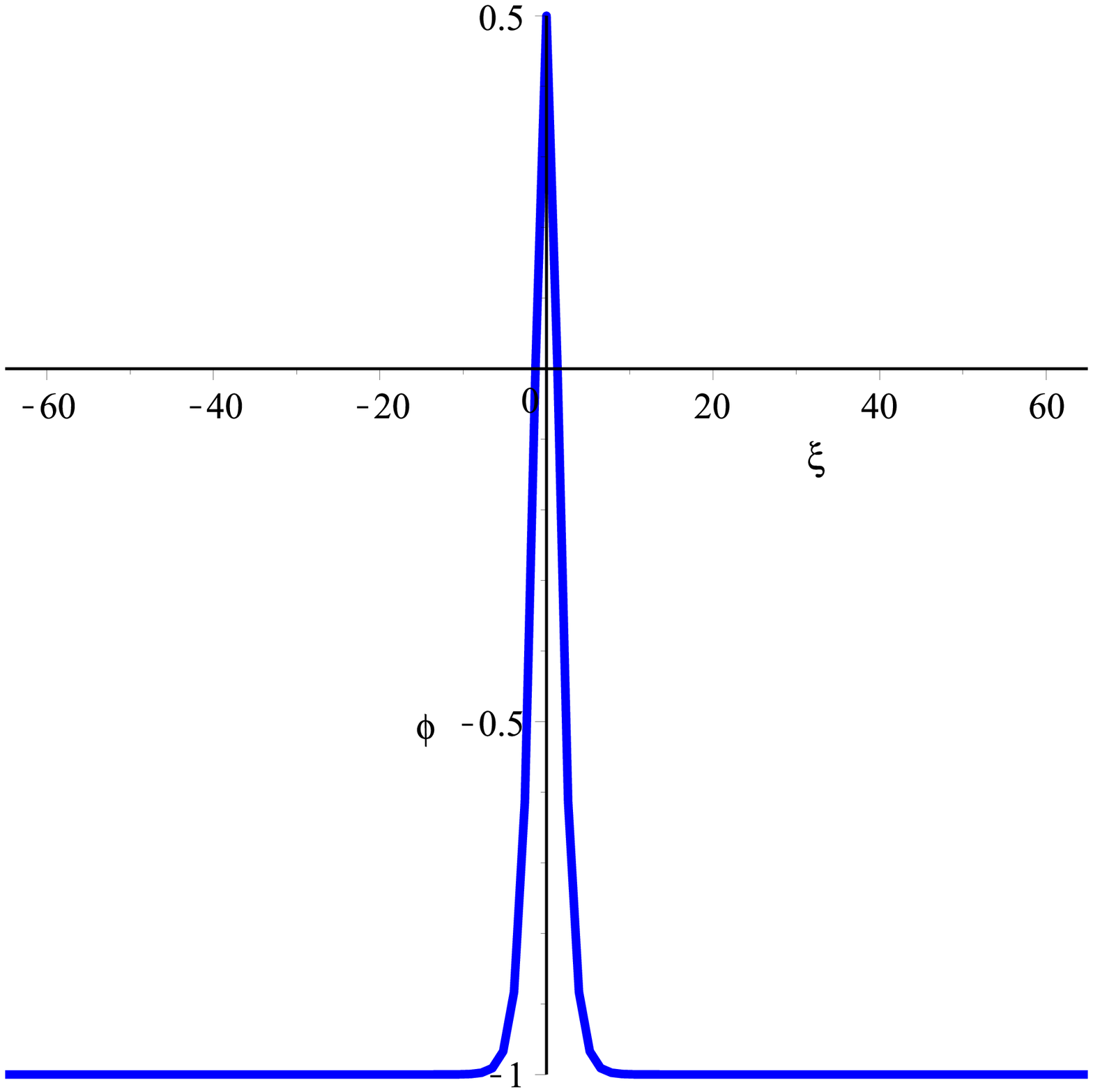}&
\epsfxsize=6.5cm \epsfysize=6.5cm \epsffile{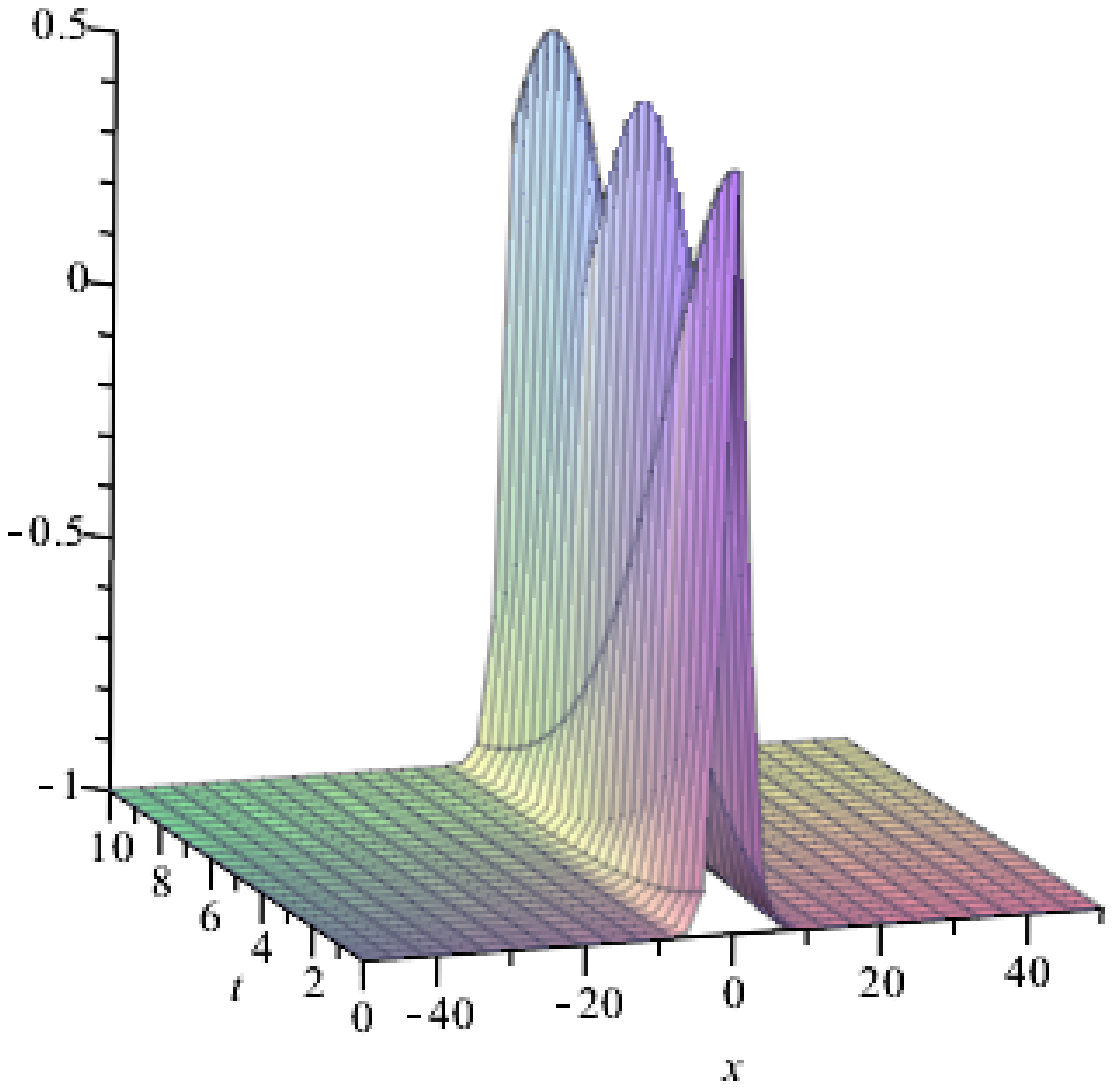}\\
\footnotesize{(a) Solitary wave  } & \footnotesize{ (b)  Bright soliton }
\end{tabular}
\end{center}
\begin{center}
\footnotesize {{Fig. 2} \ \ Solitary wave solution and bright soliton when $(c-1)^2>-2g$. }
\end{center}
Through the above analysis, we give the following theorem:
\begin{mythm}
For $(c-1)^2>-2g$, Eq. (\ref{PRLW})$|_{\tau=0}$ has a solitary wave solution given by (\ref{homoclinic2}).
\end{mythm}
\section{Persistence of solitary wave solution {with two kinds of small perturbations}}
In this section, the GSP approach is employed to analyze the persistence of solitary wave solution of Eq. (\ref{PRLW}) with KS and ME perturbations, respectively.
\subsection{{Solitary wave solution persists with KS perturbation}}
Substituting (\ref{transformations1}) into Eq. (\ref{PRLW}) with KS perturbation and integrating both sides once with respect to $\xi$,  we have
\begin{equation}\label{ODEs2aa}
(1-c)\phi+F+c\phi''+\tau(\phi'+\phi''')=g,
\end{equation}
where $F=\displaystyle\int^{\xi}_{-\infty}\psi\phi'ds$ and $\psi(\xi)=\displaystyle\int^{\infty}_0\frac{4t}{\tau^2}e^{-{\frac{2t}{\tau}}}\phi(\xi+ct)dt$.

Next, we obtain
\begin{equation}\label{function3}
\begin{array}{rcl}
\displaystyle\frac{d\psi}{d\xi}&=&\displaystyle\int^{\infty}_{0}\frac{4t}{\tau^2}e^{-{\frac{2t}{\tau}}}\phi_\xi(\xi+ct)dt\\
\vspace{.25cm}
&=&\displaystyle\frac{1}{c}\displaystyle\int^{\infty}_{0}\frac{4t}{\tau^2}e^{-{\frac{2t}{\tau}}}\phi_t(\xi+ct)dt\\
\vspace{.15cm}
&=&\displaystyle\frac{1}{c}\displaystyle\int^{\infty}_{0}\frac{4t}{\tau^2}e^{-{\frac{2t}{\tau}}}d\phi \\
\vspace{.15cm}
&=&\displaystyle\frac{1}{c}\big[\displaystyle\frac{2}{\tau}\int^{\infty}_{0}\frac{4t}{\tau^2}e^{-{\frac{2t}{\tau}}}\phi(\xi+ct)dt
-\frac{1}{\tau}\displaystyle\int^{\infty}_{0}\frac{4}{\tau}e^{-{\frac{2t}{\tau}}}\phi(\xi+ct)dt \big]\\
\vspace{.15cm}
&=&\displaystyle\frac{1}{c\tau}(2\psi-\zeta),
\end{array}
\end{equation}
where $\zeta=\displaystyle\int^{\infty}_{0}\frac{4}{\tau}e^{-{\frac{2t}{\tau}}}\phi(\xi+ct)dt$.

Similarly, differentiating both sides with respect to $\xi$, it gives
\begin{equation}\label{function4}
\frac{d\zeta}{d\xi}=\frac{2}{c\tau}(\zeta-2\phi).
\end{equation}
Thus, combining Eq. (\ref{ODEs2aa}), (\ref{function3}) and (\ref{function4}), we obtain a five-dimensional slow system as follows
\begin{equation}\label{Sys11}
\begin{cases}
\phi'=y,\\
y'=z,\\
\tau z'=(c-1)\phi-F-cz+g-\tau y,\\
c\tau\psi'=2\psi-\zeta,\\
c\tau\zeta'=2(\zeta-2\phi).
\end{cases}
\end{equation}
Setting $\tau=0$, the critical manifold is
$$M_0=\{(\phi,y,z,\psi,\zeta)\in \mathbb{R}^5\big|\psi=\frac{1}{2}\zeta=\phi,(c-1)\phi-\frac{1}{2}\phi^{2}-cz+g=0\}.$$
Taking $\chi=\tau\xi$, the corresponding fast system is
\begin{equation}\label{Sys22}
\begin{cases}
\dot{\phi}=\tau y,\\
\dot{y}=\tau z,\\
\dot{z}=(c-1)\phi-F-cz+g-\tau y,\\
\dot{\psi}=\displaystyle\frac{1}{c}(2\psi-\zeta),\\
\dot{\zeta}=\displaystyle\frac{2}{c}(\zeta-2\phi),
\end{cases}
\end{equation}
where $\cdot=\frac{d}{d\chi}$. Let $A$ be the linearized matrix of the fast system (\ref{Sys22}) restricted to $M_0$. It is of the form
\begin{equation}
A=\left(
\begin{array}{ccccc}
0 & 0 & 0 & 0 & 0\\
0 & 0 & 0 & 0 & 0 \\
\ast & 0 & -c & 0 & 0\\
0 & 0 & 0 & \frac{2}{c} & -\frac{1}{c} \\
-\frac{4}{c} & 0 & 0 & 0 & \frac{2}{c}
\end{array}
\right)
\end{equation}
that the eigenvalues of A are $0$, $0$, $-c$, $\frac{2}{c}$ and $\frac{2}{c}$. Hence, $M_0$ is normally hyperbolic manifold (see \cite{FN-JDE}). And there exists a submanifold $M_\tau$ which is a differentiable homeomorphism and locates near $M_0$ at a distance $O(\tau)$ (see \textbf{Lemma} \ref{Lemma1.1}). The invariant submanifold $M_\tau$ is represented by
\begin{equation}\label{ODEs2}
\begin{array}{rcl}
M_\tau=\{(\phi,y,z,\psi,\zeta)\in \mathbb{R}^5\big|\psi=\phi+p(\phi,y,\tau),\zeta=2\phi+q(\phi,y,\tau),\\
(c-1)\phi-\frac{1}{2}\phi^{2}-cz+g=\omega(\phi,y,\tau)\},
\end{array}
\end{equation}
where $g$, $h$ and $\omega$ are smooth enough with respect to $\tau$ and satisfy
$$p(\phi,y,0)=q(\phi,y,0)=\omega(\phi,y,0)=0.$$
Therefore, $g$, $h$ and $\omega$ can be Taylor's expansion on $\tau$
\begin{equation}\label{expansion1}
\begin{array}{rcl}
p(\phi,y,\tau)&=& \tau p_1+\tau^2 p_2+\cdots,\\
q(\phi,y,\tau)&=&\tau q_1+\tau^2 q_2+\cdots,\\
\omega(\phi,y,\tau)&=&\tau \omega_1+\tau^2 \omega_2+\cdots.\\
\end{array}
\end{equation}
Substituting (\ref{expansion1}) into (\ref{Sys11}) yields
\begin{equation}\label{expansion1a}
\begin{array}{rcl}
\tau \frac{1}{c}\big[(c-1)y-\phi y-\tau(\frac{\partial \omega_1}{\partial \phi}y+\frac{\partial \omega_1}{\partial y}z)+ O(\tau^2) \big]&=& -F+\frac{1}{2}\phi^{2}+\tau\omega_1-\tau y+O(\tau^2),\\
\tau \big[y+\tau (\frac{\partial p_1}{\partial \phi}y+\frac{\partial p_1}{\partial y}z)+ O(\tau^2)\big] &=&\tau \frac{2}{c}p_1+ O(\tau^2),\\
\tau \big[2y+\tau (\frac{\partial q_1}{\partial \phi}y+\frac{\partial q_1}{\partial y}z)+ O(\tau^2)\big] &=&\tau \frac{2}{c}q_1+ O(\tau^2).
\end{array}
\end{equation}
Comparing the coefficients of $\tau$, we have
\begin{equation}\label{expansion1b}
\begin{array}{rcl}
p_1 &=&\displaystyle\frac{c}{2}y,\\
q_1 &=&cy,\\
\omega_1 &=&\displaystyle \frac{c-1}{c}y-\frac{1}{c}\phi y+y.
\end{array}
\end{equation}
Hence, slow system restricted to $M_\tau$ is
\begin{equation}\label{critical manifold1}
\begin{cases}
\phi'=y,\\
y'=\displaystyle\frac{1}{c}\big[(c-1)\phi-\frac{1}{2}\phi^{2}+g\big]-\tau\frac{1}{c}\big(\frac{2c-1}{c}y-\frac{1}{c}\phi y \big)+O(\tau^2).
\end{cases}
\end{equation}
Secondly, let us consider the existence of solitary wave solutions. Here, Melnikov method (see \cite{MVK-TMM,GJ-NY}) is employed to detect the existence of solitary wave solutions under KS perturbation.

The Melnikov function of system (\ref{critical manifold1}) is defined by
\begin{equation}\label{Melnikov}
\begin{array}{rcl}
M_{\textmd{hom-KS}}(c,g)=\displaystyle\frac{1}{c}\displaystyle\oint_{\Gamma}(\frac{2c-1}{c}y^2+\frac{c-1}{c}\phi y^2)d\xi=\frac{2c-1}{c^2}I_1-\frac{1}{c^2}I_2,
\end{array}
\end{equation}
where $I_1=\displaystyle\oint_{\Gamma}y^2d\xi$ and $I_2=\displaystyle\oint_{\Gamma}\phi y^2d\xi$. Due to $I_1$ and $I_2$ are integrated over a closed curve, and the time variable $\xi$ is represented by the state variable $\phi$ on homoclinic orbit $\Gamma$. Thus, by (\ref{expression1}), we have
\begin{equation}\label{I1}
\begin{array}{rcl}
I_1&=&\displaystyle\oint_{\Gamma}y^2d\xi=\displaystyle\oint_{\Gamma}yd\phi=2\sqrt{\frac{1}{3c}} \displaystyle\int_{\phi_r}^{\phi_1}\sqrt{(\phi_r-\phi)(\phi-\phi_1)^2}d\phi\\
&=& 2\sqrt{\frac{1}{3c}} \displaystyle\int_{\phi_r}^{\phi_1}(\phi-\phi_1)\sqrt{\phi_r-\phi}d\phi\\
&=& -\frac{8}{45}\sqrt{\frac{3}{c}}(\phi_r-\phi_1)^{\frac{5}{2}}\\
&=&-\frac{24}{5}\sqrt{\frac{1}{c}}(c^2-2c+2g+1)^{\frac{5}{4}},
\end{array}
\end{equation}
and
\begin{equation}\label{I2}
\begin{array}{rcl}
I_2&=&\displaystyle\oint_{\Gamma}\phi y^2d\xi=\displaystyle\oint_{\Gamma}\phi yd\phi\\
&=&2\sqrt{\frac{1}{3c}} \displaystyle\int_{\phi_r}^{\phi_1}\phi\sqrt{(\phi_r-\phi)(\phi-\phi_1)^2}d\phi\\
&=& 2\sqrt{\frac{1}{3c}} \displaystyle\int_{\phi_r}^{\phi_1}\phi(\phi-\phi_1)\sqrt{\phi_r-\phi}d\phi\\
&=&2\sqrt{\frac{1}{3c}} (\displaystyle\int_{\phi_r}^{\phi_1}\phi^2\sqrt{\phi_r-\phi}d\phi-\phi_1\displaystyle\int_{\phi_r}^{\phi_1}\phi\sqrt{\phi_r-\phi}d\phi)\\
&=& -\frac{8\sqrt{3}}{315}\sqrt{\frac{1}{c}}(\phi_r-\phi_1)^{\frac{3}{2}}(4{\phi_r}^2-\phi_r\phi_1-3{\phi_1}^2)\\
&=&-\frac{8}{35}\sqrt{\frac{1}{c}}(15\sqrt{c^2-2c+2g+1}+21c-21)(c^2-2c+2g+1)^{\frac{5}{4}}.
\end{array}
\end{equation}
Therefore, we substitute (\ref{I1}) and (\ref{I2}) into (\ref{Melnikov}) to obtain the following expression
\begin{equation}\label{Hom-Melnikov1}
\begin{array}{rcl}
M_{\textmd{hom-KS}}(c,g)&=&\frac{8}{35c^2}\sqrt{\frac{1}{c}}(15\sqrt{c^2-2c+2g+1}+21c-21)(c^2-2c+2g+1)^{\frac{5}{4}}\\
&&-\frac{24(2c-1)}{5c^2}\sqrt{\frac{1}{c}}(c^2-2c+2g+1)^{\frac{5}{4}}.
\end{array}
\end{equation}
Since the complexity of expression of $M_{\textmd{hom-KS}}(c,g)$, it is hard to solve $c$ from $M_{\textmd{hom-KS}}(c,g)$ directly. The expression of $M_{\textmd{hom-KS}}(c,g)$
is simplified by
\begin{equation}\label{S-Hom-Melnikov1}
\begin{array}{rcl}
M_{\textmd{hom-KS}}(c,g)=\frac{24}{5c^2}\sqrt{\frac{1}{c}}(c^2-2c+2g+1)^{\frac{5}{4}}M^*_{\textmd{hom-KS}}(c,g),
\end{array}
\end{equation}
where $M^*_{\textmd{hom-KS}}(c,g)=\frac{5}{7}\sqrt{c^2-2c+2g+1}-c$. Because $c>0$ and $(c-1)^2>-2g$, it is easy to know that whether $M_{\textmd{hom-KS}}(c,g)=0$ has a simple zero if only if $M^*_{\textmd{hom-KS}}(c,g)=0$ has one.

Next, let us discuss the zero of $M^*_{\textmd{hom-KS}}(c,g)=0$. We have
\begin{equation}\label{D-Hom-Melnikov1}
\begin{array}{rcl}
\frac{\partial M^*_{\textmd{hom-KS}}(c,g)}{\partial c}=\frac{5(c-1)}{7\sqrt{c^2-2c+2g+1}}-1.
\end{array}
\end{equation}
To proceed, we divide into two cases.

\noindent {\bf Case} { (a) $g<0$. Since $(c-1)^2>-2g$, it implies that $c>1+\sqrt{-2g}$ or $0<c<1-\sqrt{-2g}$. For $c>1+\sqrt{-2g}$, we see that $\frac{\partial M^*_{\textmd{hom-KS}}(c,g)}{\partial c}>0$ when $c\in(1+\sqrt{-2g},1+\frac{7}{6}\sqrt{-3g})$ and $\frac{\partial M^*_{\textmd{hom-KS}}(c,g)}{\partial c}<0$ when $c\in(1+\frac{7}{6}\sqrt{-3g},+\infty)$. Obviously, $M^*_{\textmd{hom-KS}}(1+\frac{7}{6}\sqrt{-3g},g)$ is the maximum for $c>1+\sqrt{-2g}$. { Consequently,  $M^*_{\textmd{hom-KS}}(c,g)\leq \max M^*_{\textmd{hom-KS}}(c,g)=M^*_{\textmd{hom-KS}}(1+\frac{7}{6}\sqrt{-3g},g)=-\frac{4\sqrt{3}}{7}\sqrt{-g}-1<0$. Therefore, for any $c>1+\sqrt{-2g}$, the nonexistence of the simple zero of $M^*_{\textmd{hom-KS}}(c,g)$ implies  that $M_{\textmd{hom-KS}}(c,g)$ has no simple zero.
 %when $g<0$.
}

For $0<c<1-\sqrt{-2g}$ (\emph{i.e.}, $-\frac{1}{2}<g<0$), we obtain that $\frac{\partial M^*_{\textmd{hom-KS}}(c,g)}{\partial c}<0$ when $c\in(0,1-\sqrt{-2g})$. Clearly, $M^*_{\textmd{hom-KS}}(c,g)$ is monotonically decreasing about $c\in(0,1-\sqrt{-2g})$ for any $-\frac{1}{2}<g<0$. When $c\rightarrow0$, then $M^*_{\textmd{hom-KS}}(c,g)\rightarrow\frac{5}{7}\sqrt{2g+1}>0$. When $c\rightarrow1-\sqrt{-2g}$, and $M^*_{\textmd{hom-KS}}(c,g)\rightarrow-1+\sqrt{-2g}<0$. {Further, $M^*_{\textmd{hom-KS}}(c,g)$ is continuous in $c\in(0,1-\sqrt{-2g})$ for any $-\frac{1}{2}<g<0$. Therefore, the existence of a simple zero of $M^*_{\textmd{hom-KS}}(c,g)$ implies that $M_{\textmd{hom-KS}}(c,g)$ has one simple zero.}

In summary, $M_{\textmd{hom-KS}}(c,g)$ has a simple zero for any $0<c<1-\sqrt{-2g}$ when $-\frac{1}{2}< g<0$ (see Fig. 3(a)).

\noindent {\bf Case} (b) $g\geq0$. Obviously, we have $\frac{\partial M^*_{\textmd{hom-KS}}(c,g)}{\partial c}<0$ for any $c>0$ when $g\geq0$. Thus, $M^*_{\textmd{hom-KS}}(c,g)$ is monotone about $c>0$ for any $g\geq0$. Then, the existence of a simple zero of $M^*_{\textmd{hom-KS}}(c,g)$ implies that $M_{\textmd{hom-KS}}(c,g)$ has one simple zero (see Fig. 3(b) and Fig. 3(c)).}

%Then, we obtain that $\frac{\partial M^*_{\textmd{hom-KS}}(c,g)}{\partial c}>3$ for any $c\geq1+\sqrt{-2g}$ when $g<0$. Thus, $M^*_{\textmd{hom-KS}}(c,g)$ is monotonically increasing about $c\geq1+\sqrt{-2g}$ for any $g<0$. Moreover, $M^*_{\textmd{hom-KS}}(1+\sqrt{-2g},g)=1+3\sqrt{-2g}>0$, in other words, $M^*_{\textmd{hom-KS}}(c,g)=0$ has no simple zero such that $M_{\textmd{hom-KS}}(c,g)=0$ has no simple zero for any $c\geq1+\sqrt{-2g}$ when $g<0$.

%However, for $0<c<1-\sqrt{-2g}$ (\emph{i.e.}, $-\frac{1}{2}<g<0$), it needs to judge the sign of $M^*_{\textmd{hom-KS}}(c_0,g)$, where $c_0=1-\frac{21}{52}\sqrt{-13g}$ and $\frac{\partial M^*_{\textmd{hom-KS}}(c_0,g) }{\partial c}=0$.

%(i) When $M^*_{\textmd{hom-KS}}(c_0,g)<0$ (\emph{i.e.}, $-\frac{1}{2}<g<-\frac{49}{832}$) such that $M^*_{\textmd{hom-KS}}(c,g)=0$ hasn't a simple zero (see Fig. 3(a)).

%(ii) When $M^*_{\textmd{hom-KS}}(c_0,g)=0$ (\emph{i.e.}, $g=-\frac{49}{832}$), we have that $M^*_{\textmd{hom-KS}}(c,g)=0$ admits a simple zero (see Fig. 3(b)).

%(iii) When $M^*_{\textmd{hom-KS}}(c_0,g)>0$ and $M^*_{\textmd{hom-KS}}(c_1,g)<0$ (\emph{i.e.}, $-\frac{49}{832}< g<-\frac{1}{18}$), we have that $M^*_{\textmd{hom-KS}}(c,g)=0$ admits two simple zeros (see Fig. 3(c)). But, when $M^*_{\textmd{hom-KS}}(c_0,g)>0$ and $M^*_{\textmd{hom-KS}}(c_1,g)\geq0$ (\emph{i.e.}, $-\frac{1}{18}\leq g<0$), it implies that $M^*_{\textmd{hom-KS}}(c,g)=0$ admits a simple zero (see Fig. 3(d)).
\begin{center}
\begin{tabular}{ccc}
\epsfxsize=4.5cm \epsfysize=4.5cm \epsffile{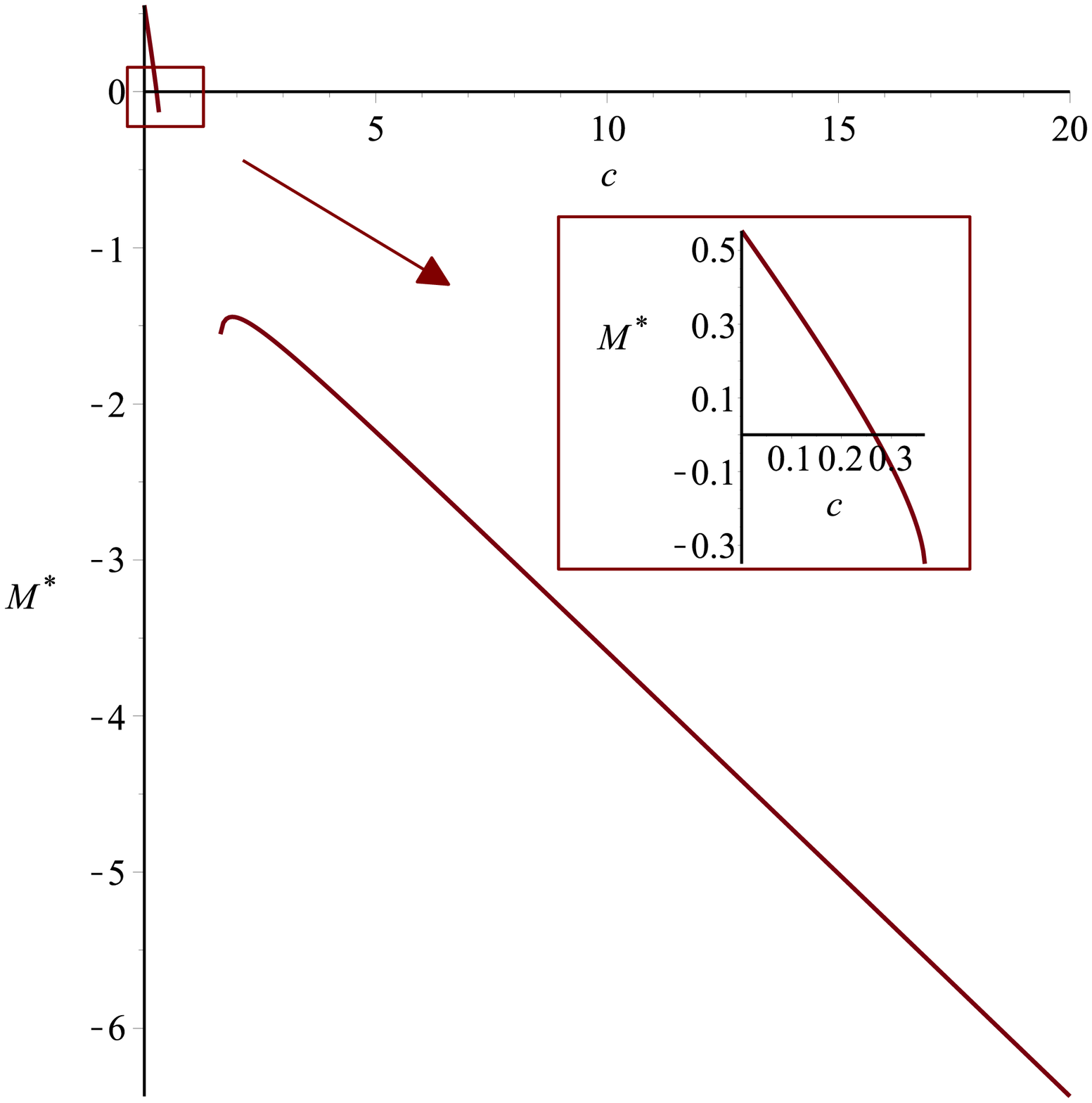}&
\epsfxsize=4.5cm \epsfysize=4.5cm \epsffile{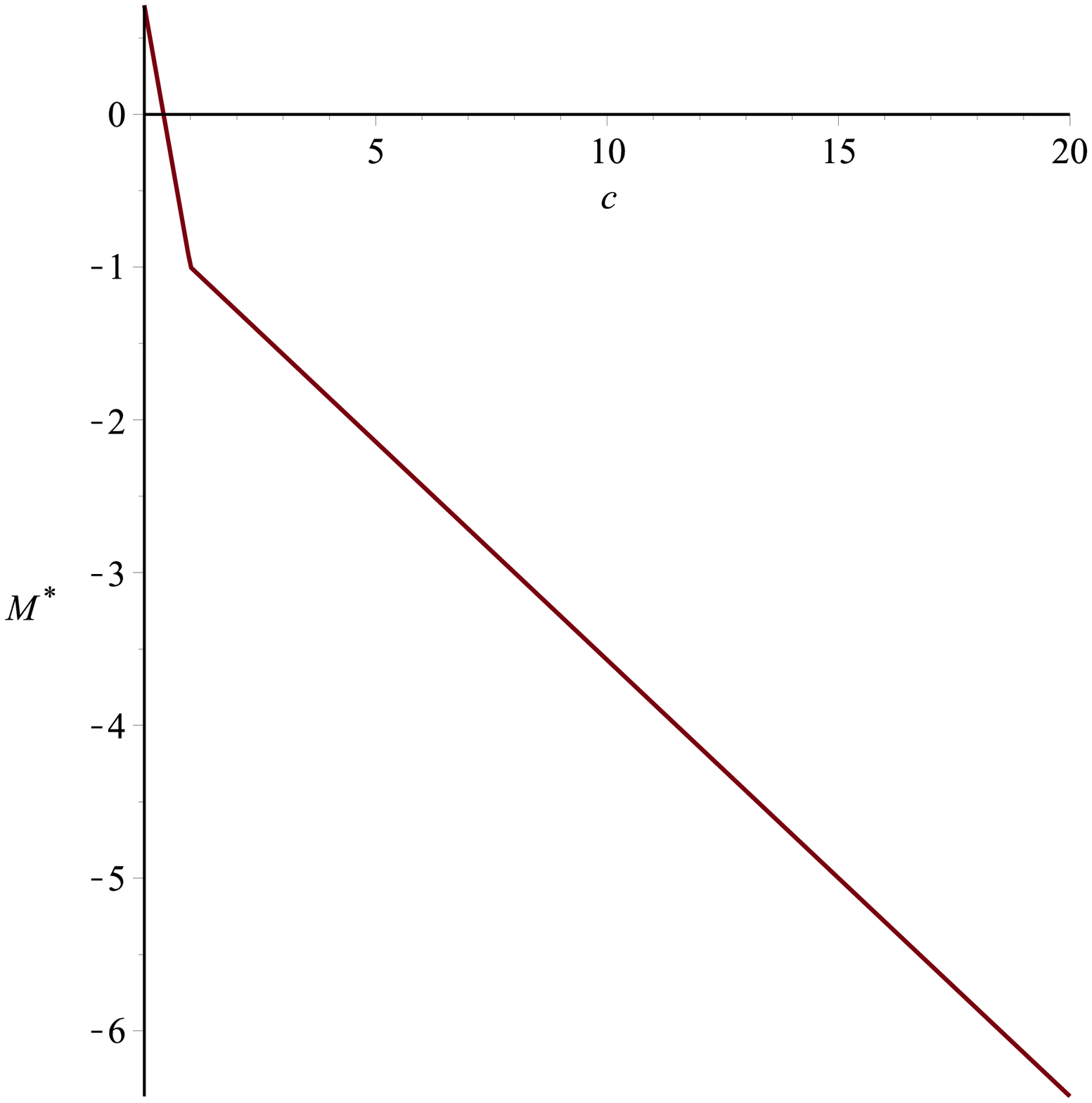}&
\epsfxsize=4.5cm \epsfysize=4.5cm \epsffile{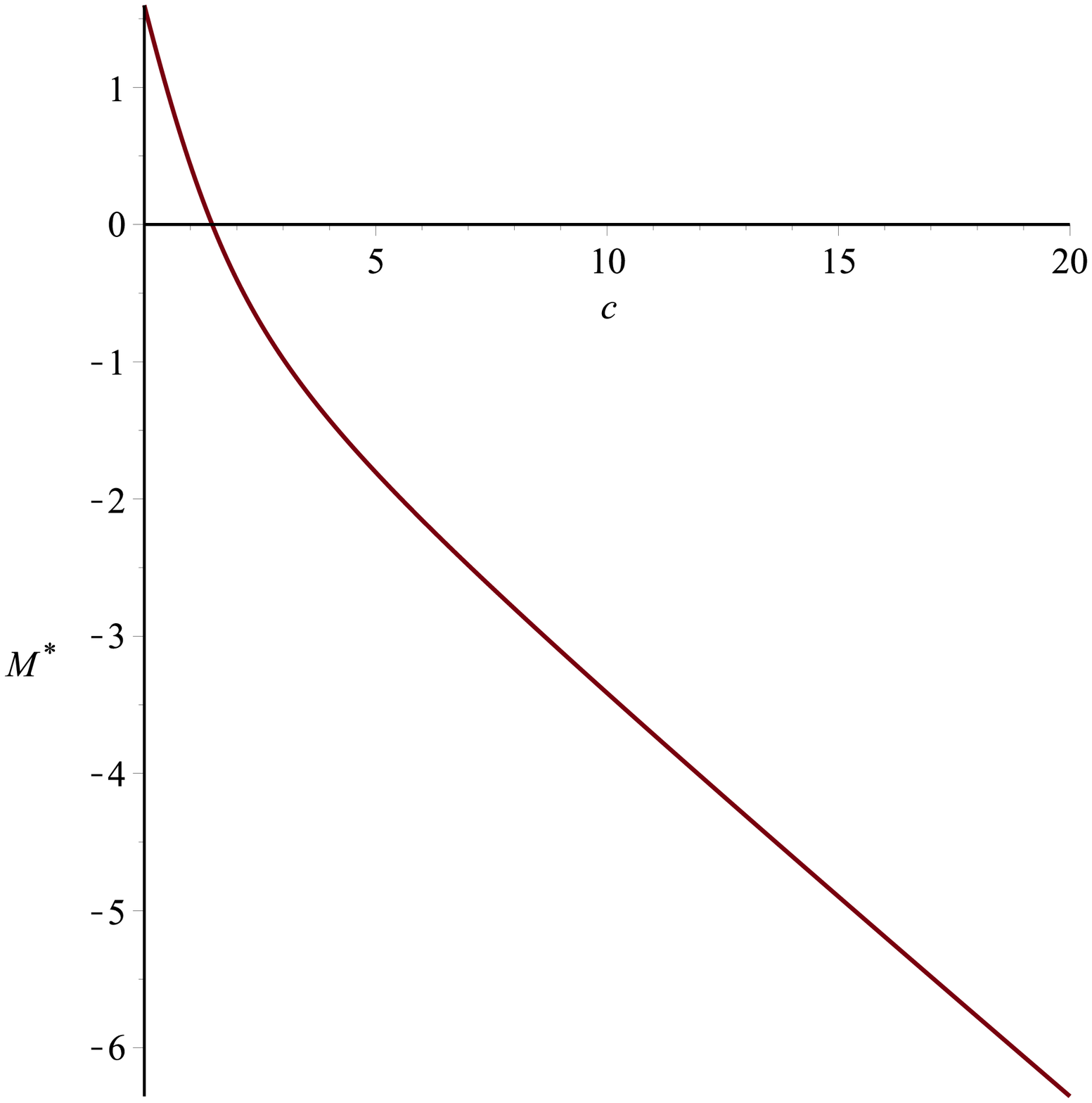}\\
\footnotesize{(a) $g=-0.2$  } & \footnotesize{ (b)  $g=0$ } &\footnotesize{(c) $g=2$  }
\end{tabular}
\end{center}
\begin{center}
\footnotesize {{Fig. 3} \ The algebraic curve of $M^*_{\textmd{hom-KS}}(c,g)$ under different value $g$.}
\end{center}
Further, we plot the algebraic curve of $M_{\textmd{hom-KS}}(c,g)$ with respect to $c$ by taking different values $g$ (see Fig. 4).
\begin{center}
\begin{tabular}{ccc}
\epsfxsize=4.5cm \epsfysize=4.5cm \epsffile{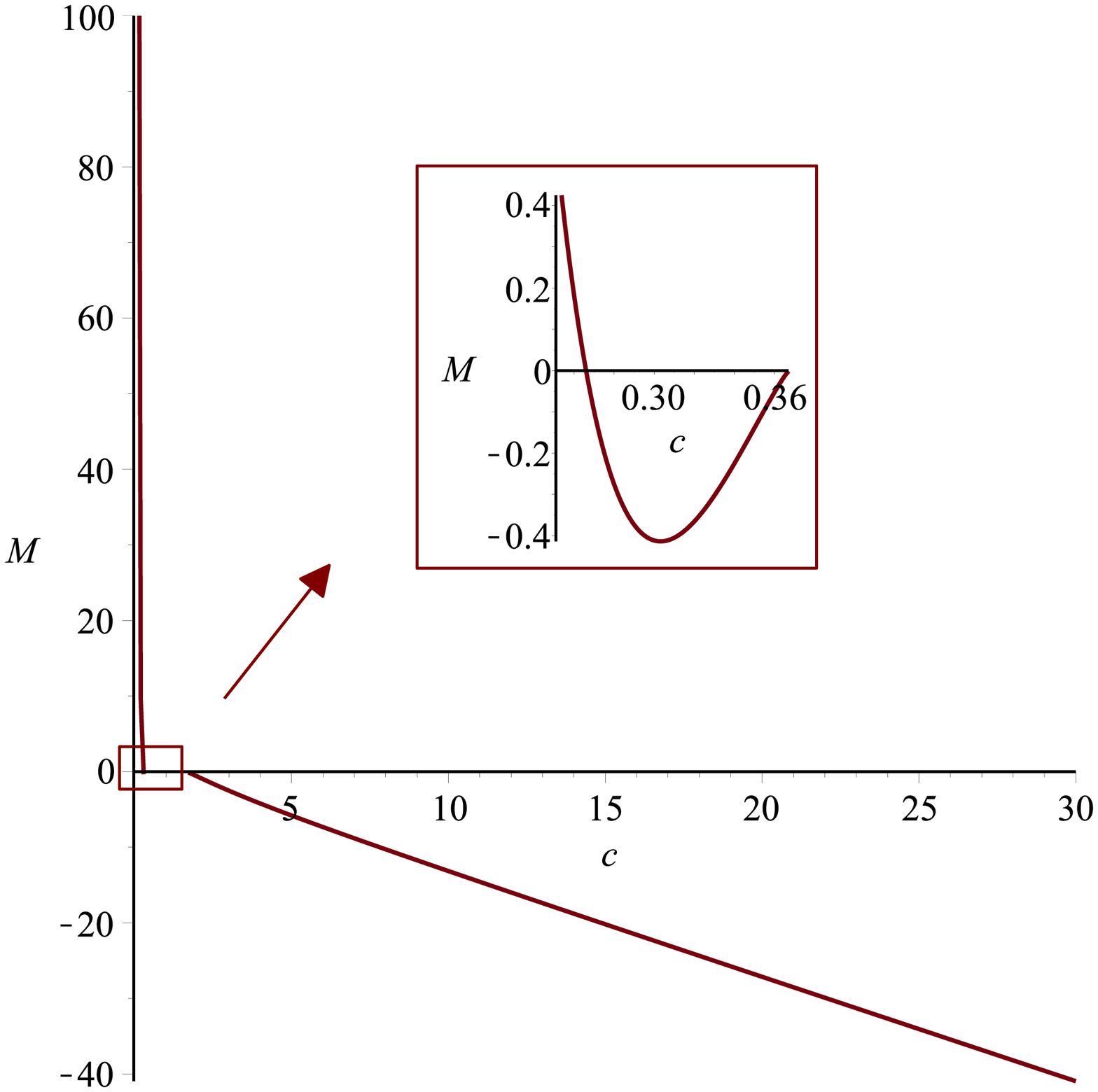}&
\epsfxsize=4.5cm \epsfysize=4.5cm \epsffile{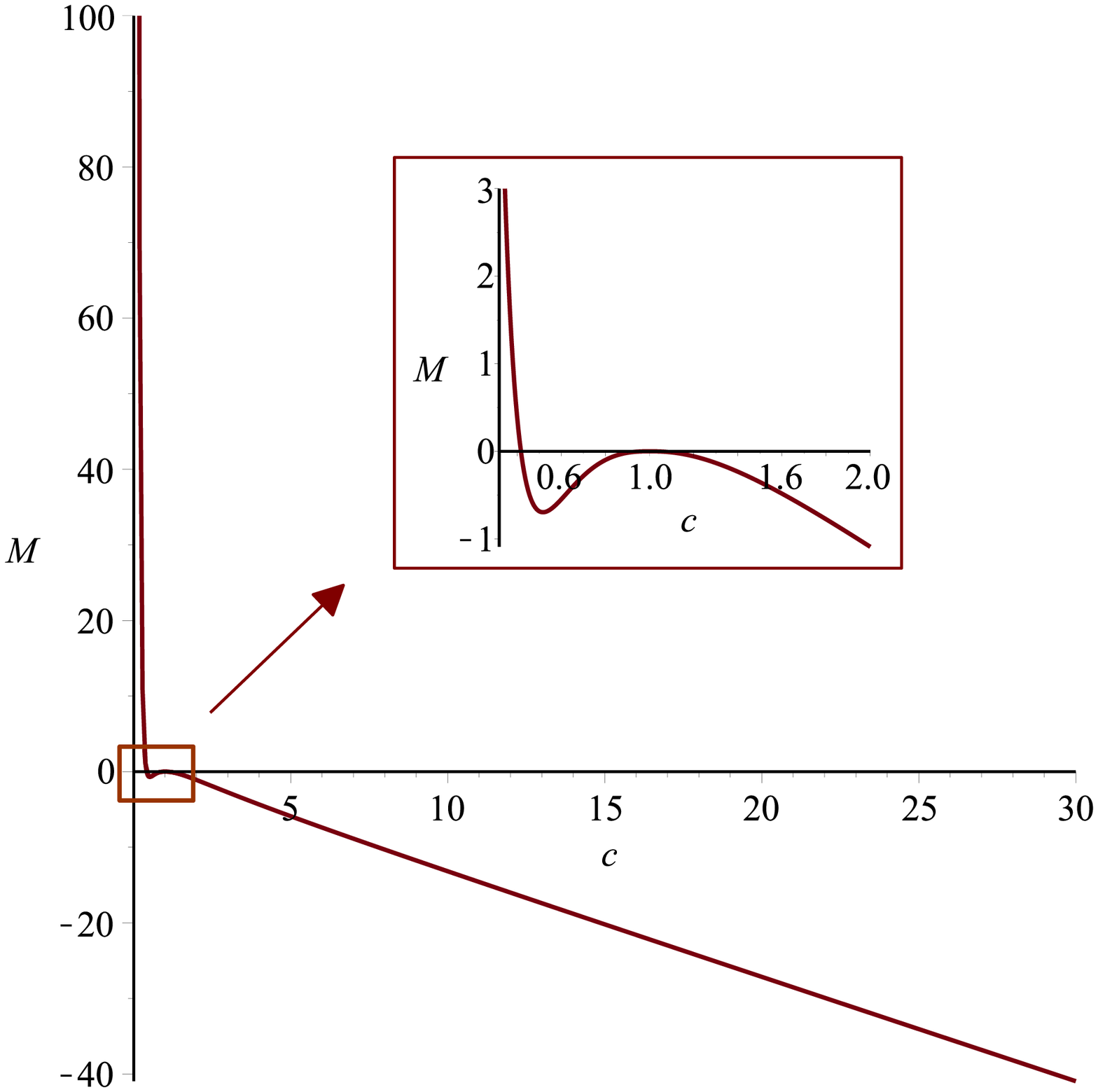}&
\epsfxsize=4.5cm \epsfysize=4.5cm \epsffile{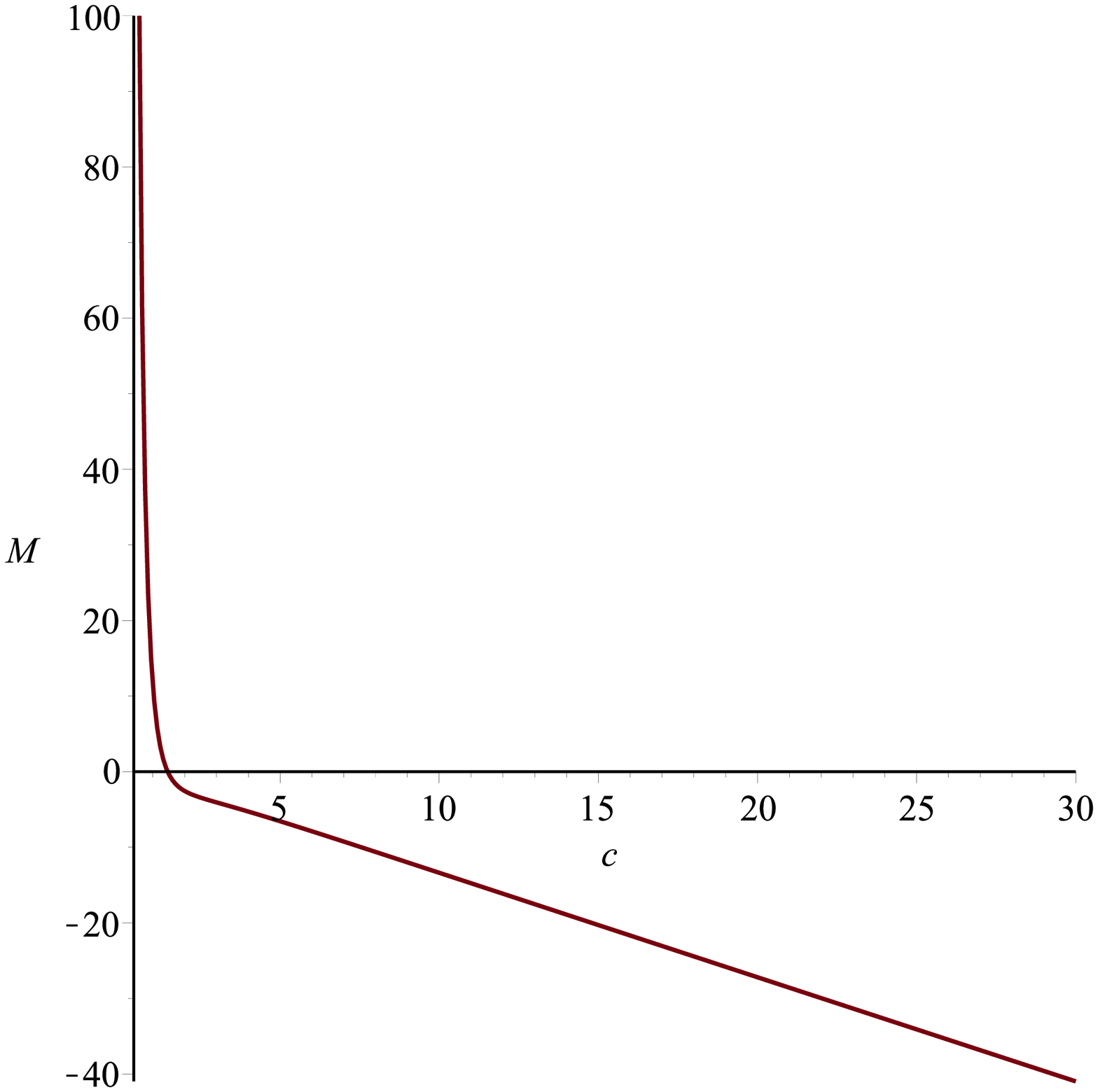}\\
\footnotesize{ (a) $g=-0.2$ } & \footnotesize{ (b) $g=0$ } & \footnotesize{ (c) $g=2$}
\end{tabular}
\end{center}
\begin{center}
\footnotesize {{Fig. 4} \ The algebraic curve of $M_{\textmd{hom-KS}}(c)$ under different value $g$.}
\end{center}
With the help of Maple, we calculate the value of $c$. The Fig. 4(b) shows that the curve intersects the $c$-axis at two points, which yields that there exist two roots $c_1=0.4166666667$ and $c_2=1$. However, if $g=0$ and $c_2=1$, it contradicts inequality $(c-1)^2>-2g$ such that we just take unique $c_1=0.4166666667$.

Therefore, there exists an unique $c^\ast>0$ such that
\begin{equation}\label{Hom-Melnikov11}
M_{\textmd{hom-KS}}(c^\ast,g)=0 \ \textmd{and} \ \frac{\partial M_{\textmd{hom-KS}}(c^\ast,g) }{\partial c}\neq0.
\end{equation}
\begin{myprop}\label{proposition2}
For any $(c-1)^2>-2g$ and $0<\tau \ll 1$, there exists a function $c(g,\tau)= c(g)+O(\tau)$ such that system
(\ref{critical manifold1}) with $c=c(g,\tau)$ possesses a homoclinic orbit near $\Gamma$.
\end{myprop}
\emph{Proof.} According to \cite[Chapter 4, Theorem 4.2.1]{Hanmaoan-book} and (\ref{Hom-Melnikov11}) that system
(\ref{critical manifold1}) with $c(g,\tau)= c(g)+O(\tau)$ has a homoclinic orbit near $\Gamma$. This completes the proof.

\begin{mythm}\label{theorem1}
For any $(c-1)^2>-2g$, there exists a function $c(g,\tau)= c(g)+O(\tau)$ such that Eq. (\ref{PRLW}) with KS perturbation has a solitary wave solution $u=u(x,t,g,\tau)$ with a wave speed $c=c(g,\tau)$. Moreover, if $\tau\rightarrow 0$, the solitary wave solution $u=u(x,t,g,\tau)$ converge to solitary wave solution (\ref{homoclinic2})  of the unperturbed Eq. $(\ref{PRLW})\mid_{\tau=0}$ with $c=c(g)$.
\end{mythm}
\subsection{{Solitary wave solution persists with ME perturbation}}
Similarly, substituting (\ref{transformations1}) into Eq. (\ref{PRLW}) with ME perturbation and integrating both sides once with respect to $\xi$,  we obtain
\begin{equation}\label{ODEs2a}
(1-c)\phi+F+c\phi''+\tau(\phi'+\phi\phi'+\phi''')=g.
\end{equation}
Therefore, combining (\ref{function3}), (\ref{function4}) and Eq. (\ref{ODEs2a}), we can also obtain a five-dimensional slow system as follows
\begin{equation}\label{Sys1}
\begin{cases}
\phi'=y,\\
y'=z,\\
\tau z'=(c-1)\phi-F-cz+g-\tau(y+\phi y),\\
c\tau\psi'=2\psi-\zeta,\\
c\tau\zeta'=2(\zeta-2\phi),
\end{cases}
\end{equation}
and it corresponds fast system is
\begin{equation}\label{Sys2}
\begin{cases}
\dot{\phi}=\tau y,\\
\dot{y}=\tau z,\\
\dot{z}=(c-1)\phi-F-cz+g-\tau(y+\phi y),\\
\dot{\psi}=\displaystyle\frac{1}{c}(2\psi-\zeta),\\
\dot{\zeta}=\displaystyle\frac{2}{c}(\zeta-2\phi).
\end{cases}
\end{equation}
Obviously, the critical manifold $M_0$ of system (\ref{Sys1}) is the same as system (\ref{Sys11}) such that it is also a normally hyperbolic manifold.

Same discussions as KS perturbation, we have
\begin{equation}\label{ODEs2}
\begin{array}{rcl}
\widetilde{M}_\tau=\{(\phi,y,z,\psi,\zeta)\in \mathbb{R}^5\big|\psi=\phi+p(\phi,y,\tau),\zeta=2\phi+q(\phi,y,\tau),\\
(c-1)\phi-\frac{1}{2}\phi^{2}-cz=\widetilde{\omega}(\phi,y,\tau)\},
\end{array}
\end{equation}
where $\widetilde{\omega}(\phi,y,\tau)=\tau  (\displaystyle \frac{c-1}{c}y-\frac{1}{c}\phi y+y+\phi y)+O(\tau^2)$.

Hence, slow system restricted to $\widetilde{M}_\tau$ is presented by
\begin{equation}\label{critical manifold}
\begin{cases}
\phi'=y,\\
y'=\displaystyle\frac{1}{c}\big[(c-1)\phi-\frac{1}{2}\phi^{2}+g\big]-\tau\frac{1}{c}\big(\frac{2c-1}{c}y+\frac{c-1}{c}\phi y \big)+O(\tau^2).
\end{cases}
\end{equation}
The Melnikov function of system (\ref{critical manifold}) is given by
\begin{equation}\label{Melnikov1}
\begin{array}{rcl}
M_{\textmd{hom-ME}}(c,g)&=&\displaystyle\frac{1}{c}\displaystyle\oint_{\Gamma}(\frac{2c-1}{c}y^2+\frac{c-1}{c}\phi y^2)d\xi=\frac{2c-1}{c^2}I_1+\frac{c-1}{c^2}I_2\\
&=&-\frac{8(c-1)}{35c^2}\sqrt{\frac{1}{c}}(15\sqrt{c^2-2c+2g+1}+21c-21)(c^2-2c+2g+1)^{\frac{5}{4}}\\
&&-\frac{24(2c-1)}{5c^2}\sqrt{\frac{1}{c}}(c^2-2c+2g+1)^{\frac{5}{4}}.
\end{array}
\end{equation}
The expression of $M_{\textmd{hom-ME}}(c,g)$
is simplified by
\begin{equation}\label{S-Hom-Melnikov1}
\begin{array}{rcl}
M_{\textmd{hom-ME}}(c,g)=-\frac{24}{5c^2}\sqrt{\frac{1}{c}}(c^2-2c+2g+1)^{\frac{5}{4}}M^*_{\textmd{hom-ME}}(c,g),
\end{array}
\end{equation}
where $M^*_{\textmd{hom-ME}}(c,g)=\frac{5(c-1)}{7}\sqrt{c^2-2c+2g+1}+c^2$. Similarly, we see that whether $M_{\textmd{hom-ME}}(c,g)=0$ has a simple zero if only if $M^*_{\textmd{hom-ME}}(c,g)=0$ has one.

Next, let us discuss the zero of $M^*_{\textmd{hom-ME}}(c,g)=0$. We have
\begin{equation}\label{D-Hom-Melnikov1-ME}
\begin{array}{rcl}
\frac{\partial M^*_{\textmd{hom-ME}}(c,g)}{\partial c}=\frac{2}{7}\frac{7c\sqrt{c^2-2c+2g+1}+5(c-1)^2+5g}{\sqrt{c^2-2c+2g+1}}.
\end{array}
\end{equation}
To proceed, we divided into two cases.

\noindent { {\bf Case} ($\widetilde{a}$) $g<0$. Since $(c-1)^2>-2g$, it implies that $c>1+\sqrt{-2g}$ or $0<c<1-\sqrt{-2g}$. For $c>1+\sqrt{-2g}$, we obtain that $\frac{\partial M^*_{\textmd{hom-ME}}(c,g)}{\partial c}>0$ when $c\in(1+\sqrt{-2g},+\infty)$. Obviously, $M^*_{\textmd{hom-ME}}(1+\sqrt{-2g},g)$ is the minimum for $c>1+\sqrt{-2g}$. { Consequently,  $M^*_{\textmd{hom-ME}}(c,g)\geq \min M^*_{\textmd{hom-ME}}(c,g)=M^*_{\textmd{hom-ME}}(1+\sqrt{-2g},g)=(1+\sqrt{-2g})^2>0$. Therefore, for any $c>1+\sqrt{-2g}$, the nonexistence of the simple zero of $M^*_{\textmd{hom-ME}}(c,g)$ implies  that $M_{\textmd{hom-ME}}(c,g)$ has no simple zero.
 %when $g<0$.
}

For $0<c<1-\sqrt{-2g}$ (\emph{i.e.}, $-\frac{1}{2}<g<0$), it obtains that $\frac{\partial M^*_{\textmd{hom-ME}}(c,g)}{\partial c}>0$ when $c\in(0,1-\sqrt{-2g})$. Clearly, $M^*_{\textmd{hom-ME}}(c,g)$ is monotonically decreasing about $0<c<1-\sqrt{-2g}$ for any $-\frac{1}{2}<g<0$. When $c\rightarrow0$, then $M^*_{\textmd{hom-ME}}(c,g)\rightarrow-\frac{5}{7}\sqrt{2g+1}<0$. When $c\rightarrow+1-\sqrt{-2g}$, and $M^*_{\textmd{hom-ME}}(c,g)\rightarrow(-1+\sqrt{-2g})^2>0$. {Further, $M^*_{\textmd{hom-ME}}(c,g)$ is continuous in $c\in(0,1-\sqrt{-2g})$ for any $-\frac{1}{2}<g<0$. Then, the existence of a simple zero of $M^*_{\textmd{hom-ME}}(c,g)$ implies that $M_{\textmd{hom-ME}}(c,g)$ has one simple zero.

In summary, $M_{\textmd{hom-ME}}(c,g)$ has a simple zero for any $0<c<1-\sqrt{-2g}$ when $-\frac{1}{2}< g<0$ (see Fig. 5(a)).}

\noindent   {\bf Case} ($\widetilde{b}$) $g\geq0$. Obviously, we have that $\frac{\partial M^*_{\textmd{hom-ME}}(c,g)}{\partial c}>0$ for any $c>0$ when $g\geq0$. Clearly, $M^*_{\textmd{hom-ME}}(c,g)$ is monotonically decreasing about $c>0$ when $g\geq0$. When $c\rightarrow0$, then $M^*_{\textmd{hom-ME}}(c,g)\rightarrow-\frac{5}{7}\sqrt{2g+1}<0$. When $c\rightarrow+\infty$, and $M^*_{\textmd{hom-ME}}(c,g)\rightarrow+\infty$. {Moreover, $M^*_{\textmd{hom-ME}}(c,g)$ is continuous in $c\in(0,+\infty)$ for any $g\geq0$. Thus, the existence of a simple zero of $M^*_{\textmd{hom-ME}}(c,g)$ implies that $M_{\textmd{hom-ME}}(c,g)$ has one simple zero.

We conclude that $M_{\textmd{hom-ME}}(c,g)$ has a simple zero for any $c>0$ when $g\geq0$ (see Fig. 5(b) and Fig. 5(c)).}
\begin{center}
\begin{tabular}{ccc}
\epsfxsize=4.5cm \epsfysize=4.5cm \epsffile{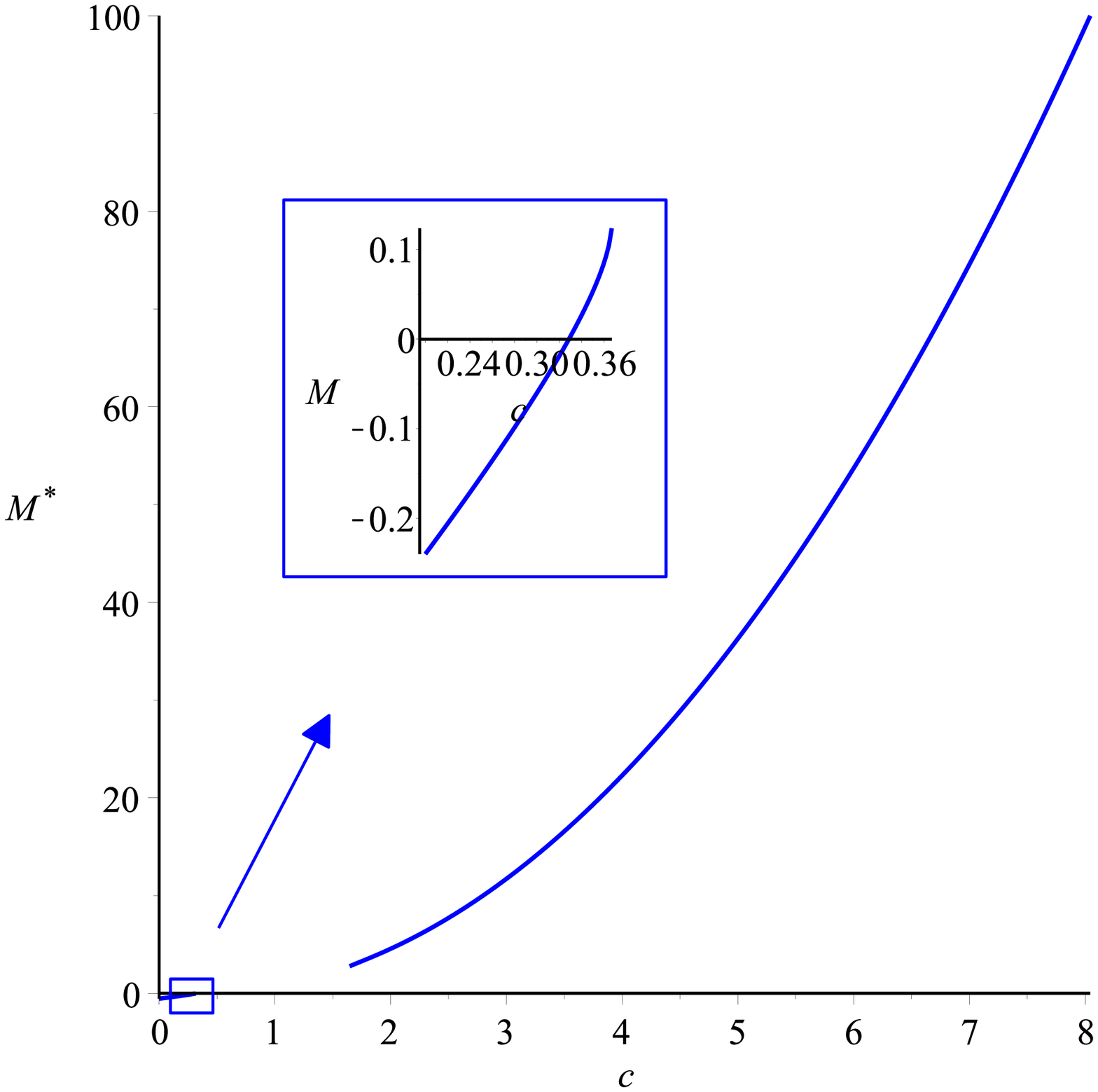}&
\epsfxsize=4.5cm \epsfysize=4.5cm \epsffile{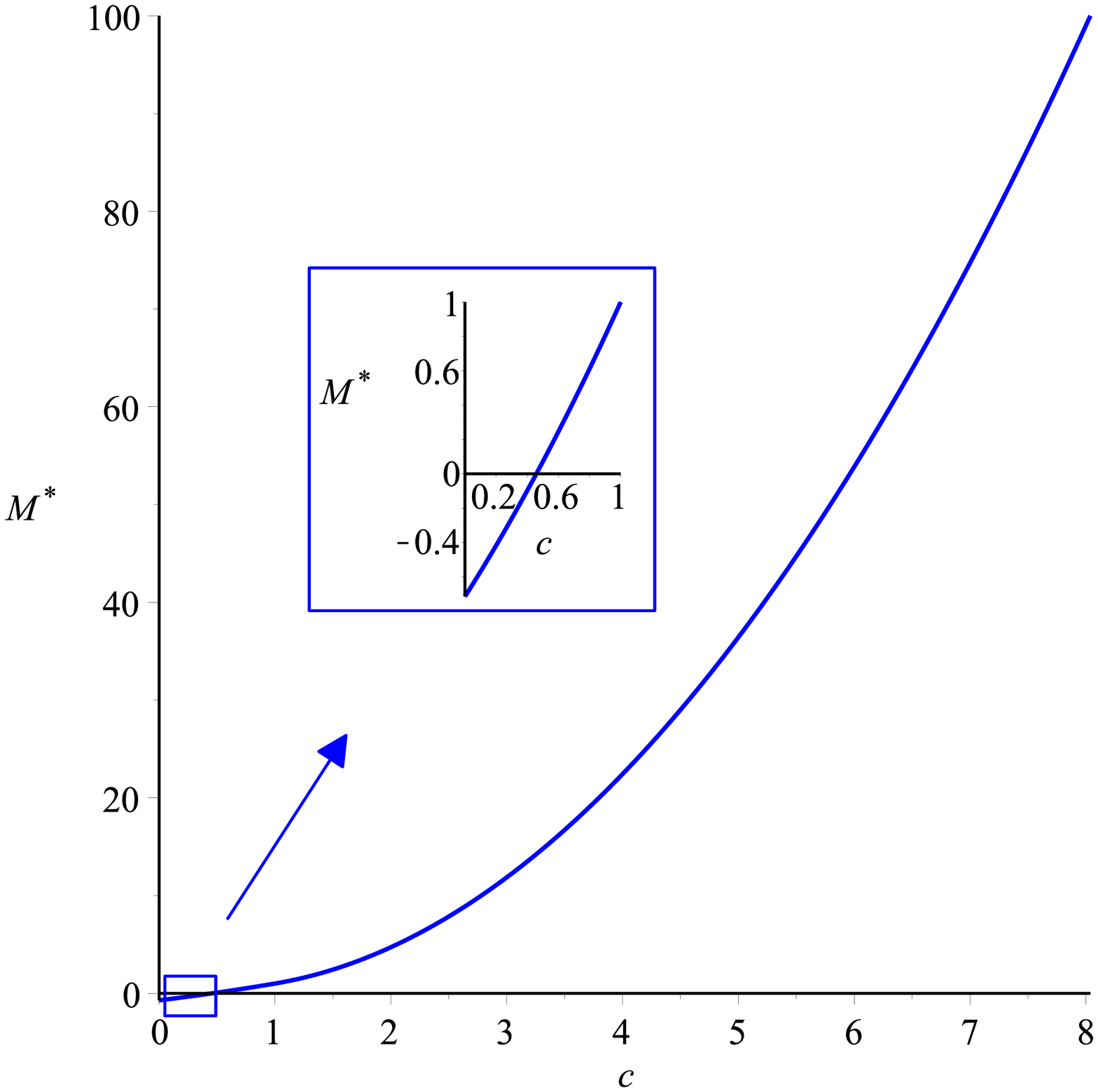}&
\epsfxsize=4.5cm \epsfysize=4.5cm \epsffile{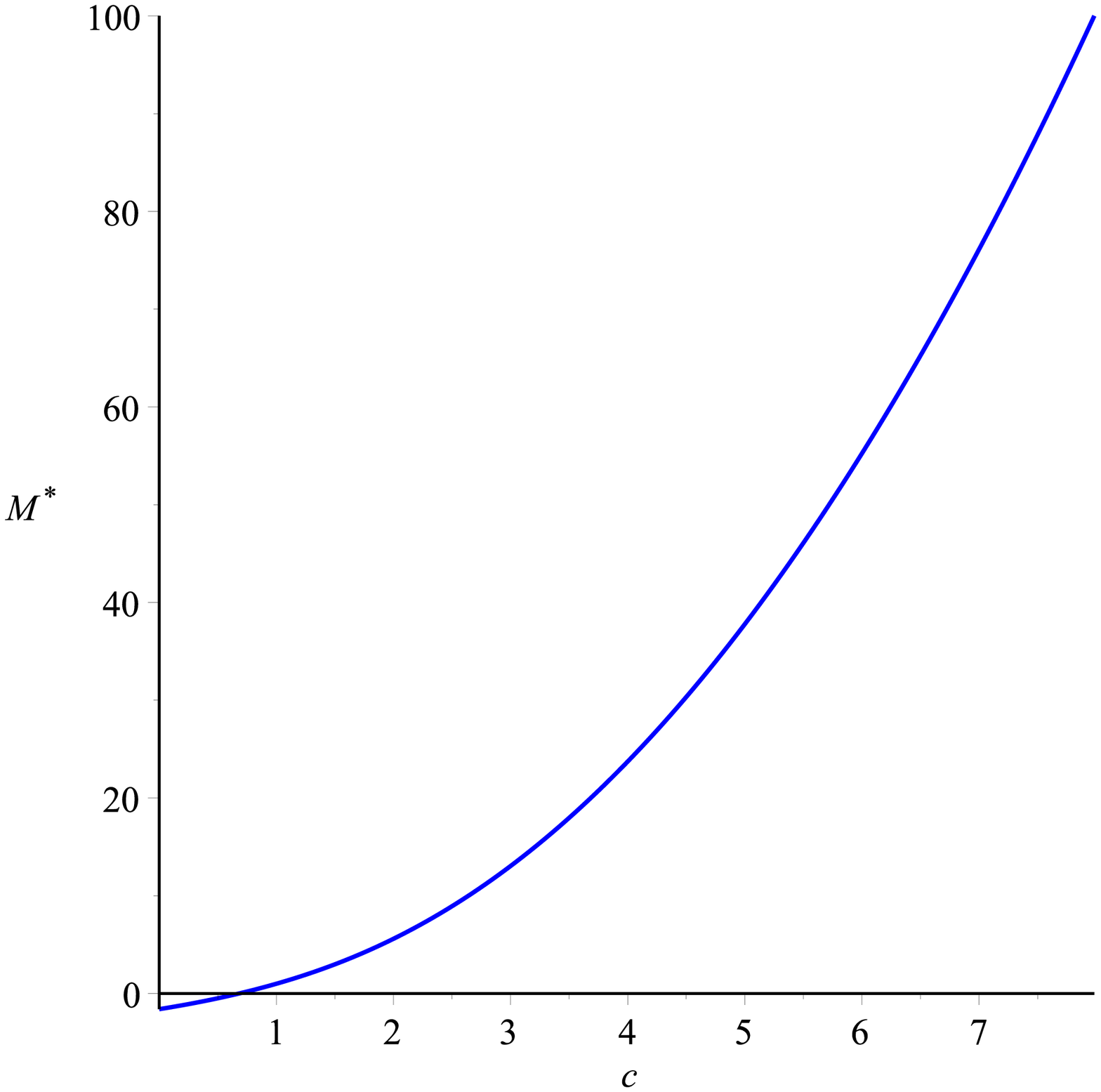}\\
\footnotesize{(a) $g=-0.2$  } & \footnotesize{ (b)  $g=0$ } &\footnotesize{(c) $g=2$  }
\end{tabular}
\end{center}
\begin{center}
\footnotesize {{Fig. 5} \ The algebraic curve of $M^*_{\textmd{hom-ME}}(c,g)$ under different value $g$.}
\end{center}
Then, we plot the algebraic curve of $M_{\textmd{hom-ME}}(c,g)$ with respect to $c$ by taking different values $g$ (see Fig. 6).
\begin{center}
\begin{tabular}{ccc}
\epsfxsize=4.5cm \epsfysize=4.5cm \epsffile{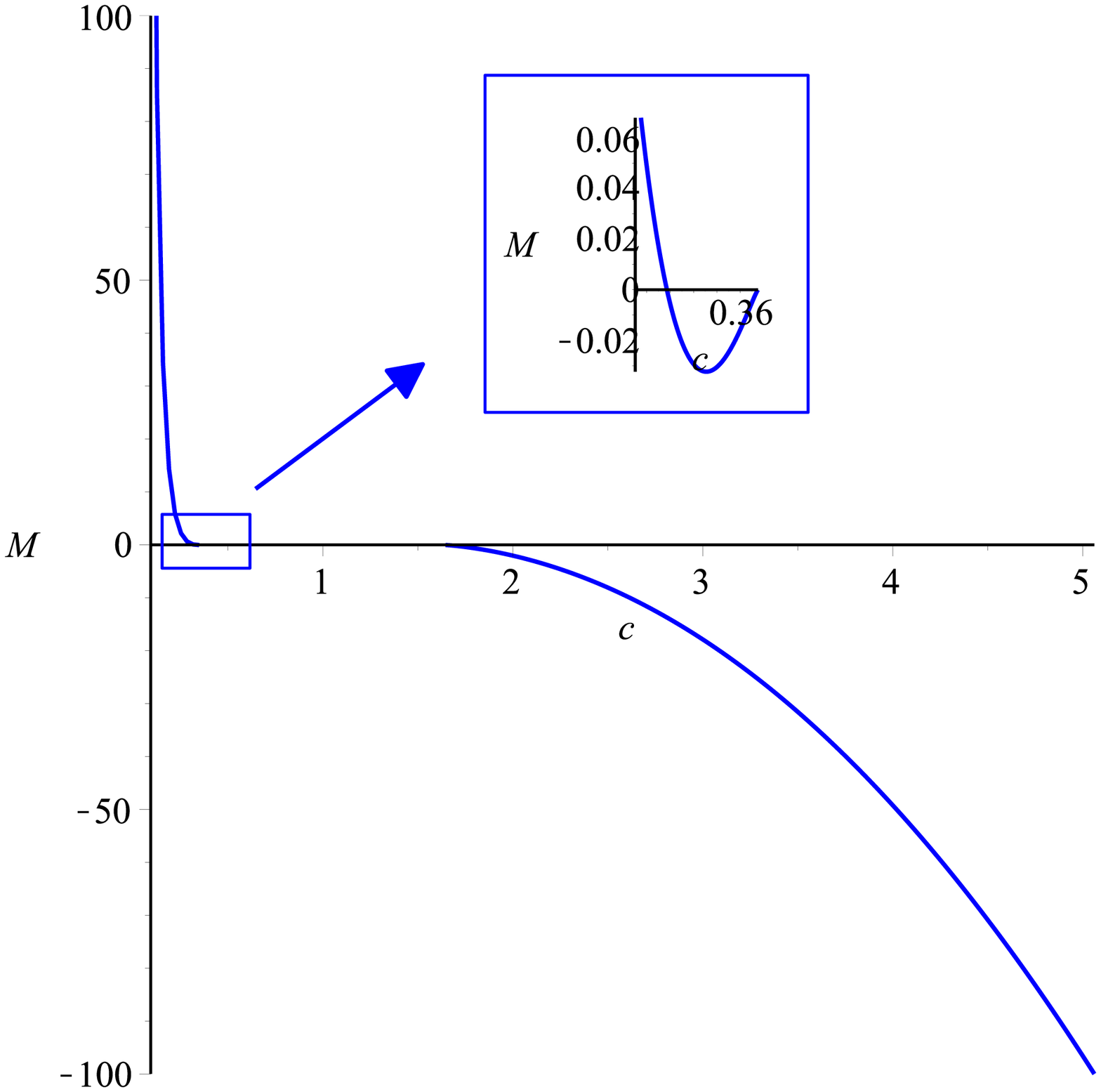}&
\epsfxsize=4.5cm \epsfysize=4.5cm \epsffile{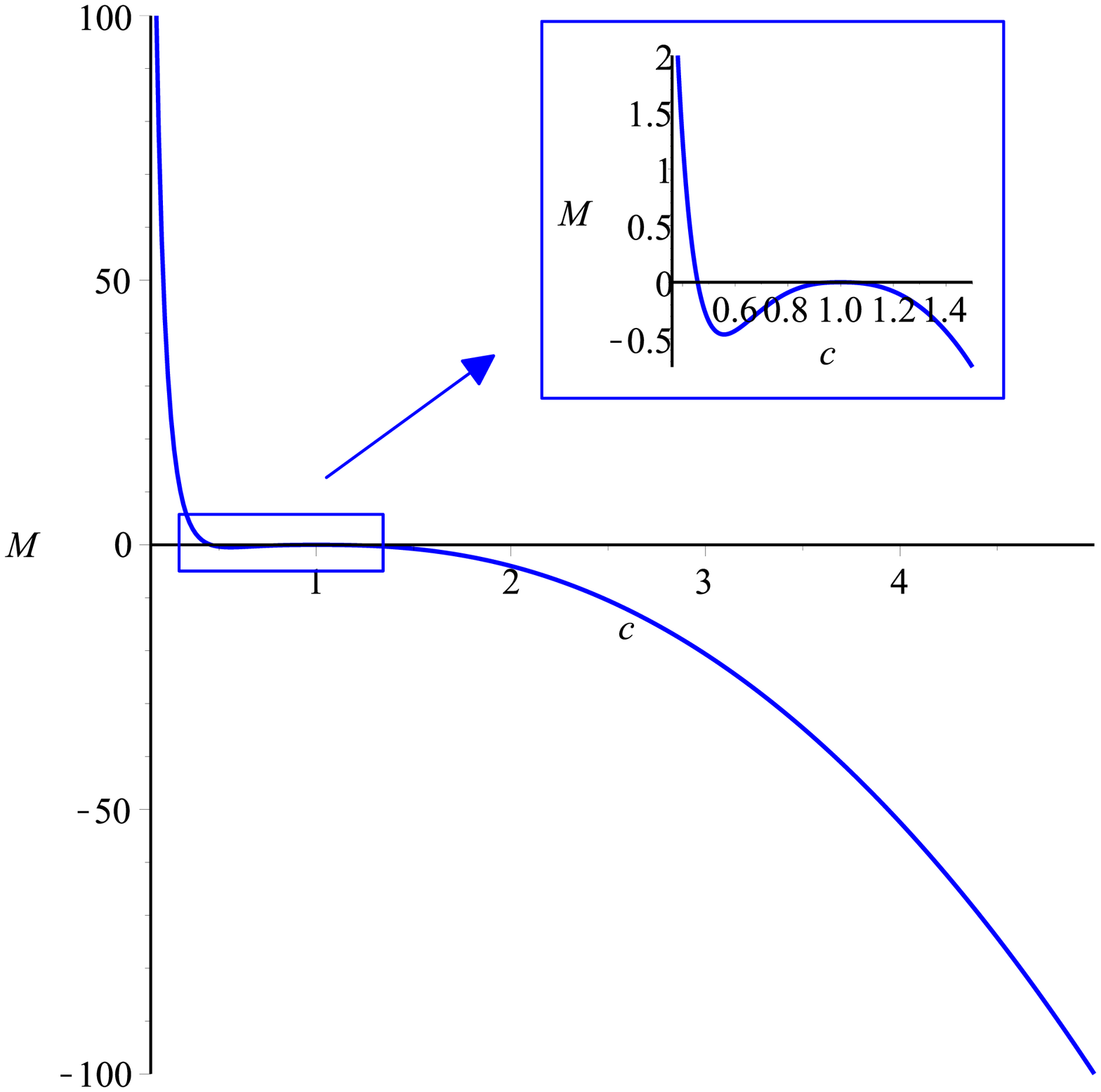}&
\epsfxsize=4.5cm \epsfysize=4.5cm \epsffile{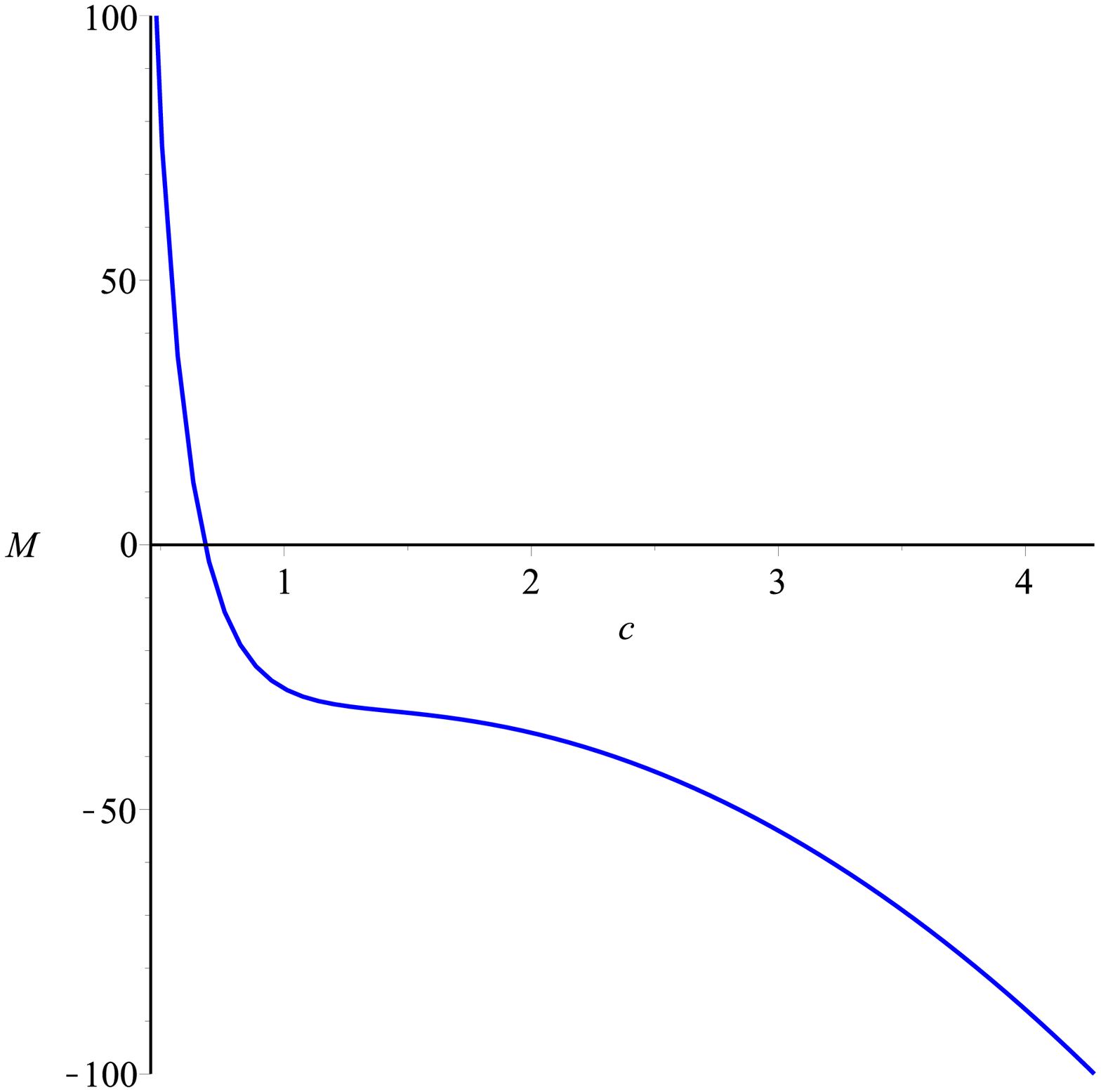}\\
\footnotesize{ (a) $g=-0.2$ } & \footnotesize{ (b) $g=0$ } & \footnotesize{ (c) $g=2$}
\end{tabular}
\end{center}
\begin{center}
\footnotesize {{Fig. 6} \ The algebraic curve of $M_{\textmd{hom-ME}}(c)$ under different value $g$.}
\end{center}
The same discussion as before, there exists an unique $c^\ast$ such that
\begin{equation}\label{Hom-Melnikova11}
M_{\textmd{hom-ME}}(c^\ast,g)=0 \ \textmd{and} \ \frac{\partial M_{\textmd{hom-ME}}(c^\ast,g) }{c}\neq0.
\end{equation}
\begin{myprop}\label{proposition2a}
For any $(c-1)^2>-2g$ and $0<\tau \ll 1$, there exists a function $c(g,\tau)= c(g)+O(\tau)$ such that system
(\ref{critical manifold}) with $c=c(g,\tau)$ possesses a homoclinic orbit near $\Gamma$.
\end{myprop}
\emph{Proof.} Similarly, according to \cite[Chapter 4, Theorem 4.2.1]{Hanmaoan-book} and (\ref{Hom-Melnikova11}) that system
(\ref{critical manifold}) with $c(g,\tau)= c(g)+O(\tau)$ has a homoclinic orbit near $\Gamma$. This completes the proof.

\begin{mythm}\label{theorem1a}
For any $(c-1)^2>-2g$, there exists a function $c(g,\tau)= c(g)+O(\tau)$ such that Eq. (\ref{PRLW}) with ME perturbation has a solitary wave solution $u=u(x,t,g,\tau)$ with a wave speed $c=c(g,\tau)$. Moreover, if $\tau\rightarrow 0$, the solitary wave solution $u=u(x,t,g,\tau)$ converge to solitary wave solution (\ref{homoclinic2}) of the unperturbed Eq. $(\ref{PRLW})\mid_{\tau=0}$ with $c=c(g)$.
\end{mythm}

\section{Numerical simulations, comparison and conclusions }
In this section, numerical simulations by maple are employed to verify the theoretical results given in previous sections. Taking $\tau=0.01$, we have

\textbf{Case I}. KS perturbation

Let $g=-0.2$, $g=0$ and $g=2$. By (\ref{Hom-Melnikov1}), we obtain $c(-0.2)\approx0.2661295689$, $c(0)\approx0.4166666667$ and $c(2)\approx1.466998871$, respectively. According to \textbf{Theorem} \ref{theorem1}, taking $c=c(-0.2)+0.0001=0.2660295689+0.0001$ and $c=c(0)+0.0001=0.4166666667+0.0001$, we set initial value to be $(\phi(0),y(0))=(\phi_r-10^{-4},0)$ which the homoclinic orbit would pass through. The phase portraits $(\phi,y)$ and time history curves $(\xi,\phi)$  of system (\ref{critical manifold1}) are plotted in Fig. (a) and Fig. (c) of table 1. Similarly, taking $c=c(2)+0.0001=1.466998871+0.0001$, we set initial value to be $(\phi(0),y(0))=(\phi_r-10^{-3},0)$, the phase portraits $(\phi,y)$ and time history curves $(\xi,\phi)$  of system (\ref{critical manifold1}) are plotted in Fig. (e) of table 1.

{The results show us that the homoclinic orbit of system (\ref{critical manifold1}) and solitary wave of Eq. (\ref{PRLW}) with suitable speed $c$ still persist under the KS perturbation.} Thus, we verify the correctness of \textbf{Proposition} \ref{proposition2} and \textbf{Theorem} \ref{theorem1}.

%\begin{center}
%\begin{tabular}{cc}
%\epsfxsize=6.5cm \epsfysize=6.5cm \epsffile{fig12a.eps}&
%\epsfxsize=6.5cm \epsfysize=6.5cm \epsffile{fig12b.eps}\\
%\footnotesize{(a)\ $c=0.2660295689+0.0001$ } & \footnotesize{(b)\ $(\xi,\phi)$ }
%\end{tabular}
%\end{center}
%\begin{center}
%\begin{tabular}{cc}
%\epsfxsize=6.5cm \epsfysize=6.5cm \epsffile{fig10a.eps}&
%\epsfxsize=6.5cm \epsfysize=6.5cm \epsffile{fig10b.eps}\\
%\footnotesize{(c)\ $c=0.4166666667+0.0001$ } & \footnotesize{(d)\ $(\xi,\phi)$ }
%\end{tabular}
%\end{center}
%\begin{center}
%\begin{tabular}{cc}
%\epsfxsize=6.5cm \epsfysize=6.5cm \epsffile{fig11a.eps}&
%\epsfxsize=6.5cm \epsfysize=6.5cm \epsffile{fig11b.eps}\\
%\footnotesize{(e)\ $c=1.466998871+0.0001$ } & \footnotesize{(f)\ $(\xi,\phi)$ }
%\end{tabular}
%\end{center}
%\begin{center}
%\footnotesize {{Fig. 7} \ Phase portraits and time history curves of KS perturbation for $\tau=0.01$. Initial value are $(\phi(0),y(0))=(\phi_r-10^{-4},0)$ in (a)-(d) and %$(\phi(0),y(0))=(\phi_r-10^{-3},0)$ in (e) and (f). }
%\end{center}

\textbf{Case II}. ME perturbation

Set $g=-0.2$, $g=0$ and $g=2$. Using (\ref{Melnikov1}), we have $c(-0.2)\approx0.3286393802$, $c(0)\approx 0.4580398915$ and $c(2)\approx0.6801960271$, respectively. By \textbf{Theorem} \ref{theorem1a}, taking  $c=c(-0.2)+0.0001=0.3286393802+0.0001$ and $c=c(0)+0.0001=0.4580398915+0.0001$, we set initial value to be $(\phi(0),y(0))=(\phi_r-10^{-4},0)$ which the homoclinic orbit would pass through. The phase portraits $(\phi,y)$ and time history curves $(\xi,\phi)$  of system (\ref{critical manifold}) are presented in in Fig. (b) and Fig. (d) of table 1. Similarly, taking $c=c(2)+0.0001=0.6801960271+0.0001$, we set initial value to be $(\phi(0),y(0))=(\phi_r-10^{-3},0)$, the phase portraits $(\phi,y)$ and time history curves $(\xi,\phi)$  of system (\ref{critical manifold}) are plotted in Fig. (f) of table 1.

Similarly, the homoclinic orbit of travelling wave system (\ref{critical manifold}) and solitary wave of Eq. (\ref{PRLW}) with suitable speed $c$ still persist under the ME perturbation. Thus, the correctness of \textbf{Proposition} \ref{proposition2a} and \textbf{Theorem} \ref{theorem1a} are verified.

%\begin{center}
%\begin{tabular}{cc}
%\epsfxsize=6.5cm \epsfysize=6.5cm \epsffile{fig7a.eps}&
%\epsfxsize=6.5cm \epsfysize=6.5cm \epsffile{fig7b.eps}\\
%\footnotesize{(a)\ $c=0.3286393802+0.0001$ } & \footnotesize{(b)\ $(\xi,\phi)$ }
%\end{tabular}
%\end{center}
%\begin{center}
%\begin{tabular}{cc}
%\epsfxsize=6.5cm \epsfysize=6.5cm \epsffile{fig5a.eps}&
%\epsfxsize=6.5cm \epsfysize=6.5cm \epsffile{fig5b.eps}\\
%\footnotesize{(c)\ $c=0.4580398915+0.0001$ } & \footnotesize{(d)\ $(\xi,\phi)$ }
%\end{tabular}
%\end{center}
%\begin{center}
%\begin{tabular}{cc}
%\epsfxsize=6.5cm \epsfysize=6.5cm \epsffile{fig6a.eps}&
%\epsfxsize=6.5cm \epsfysize=6.5cm \epsffile{fig6b.eps}\\
%\footnotesize{(e)\ $c=0.6801960271+0.0001$ } & \footnotesize{(f)\ $(\xi,\phi)$ }
%\end{tabular}
%\end{center}
%\begin{center}
%\footnotesize {{Fig. 8} \ Phase portraits and time history curves of ME perturbation for $\tau=0.01$. Initial value are $(\phi(0),y(0))=(\phi_r-10^{-4},0)$ in (a)-(d) and $(\phi(0),y(0))=(\phi_r-10^{-3},0)$ in (e) and (f). }
%\end{center}
\begin{table}[H]
\begin{center}
\caption{Phase portraits of travelling system (\ref{critical manifold1}) and time history are presented for two kinds of perturbations under $\tau=0.01$.}
\label{table:1}
\begin{tabular}{c|c|c}
\hline   \textbf{$g$} & \textbf{The perturbation of KS} & \textbf{The perturbation of ME} \\
\hline $g=-0.2$ &
\begin{tabular}{cc}
\epsfxsize=3cm \epsfysize=3cm \epsffile{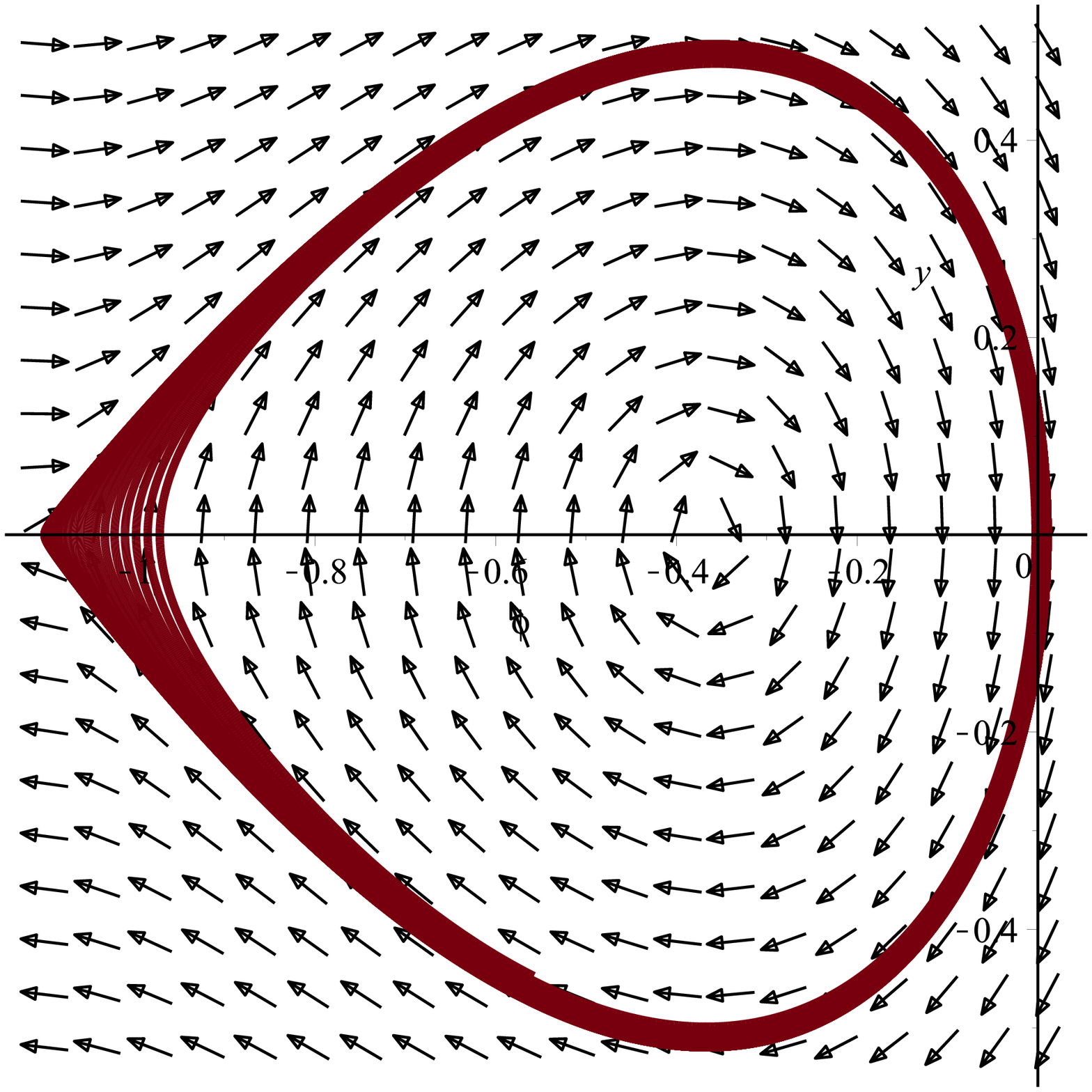}\ \ \ \epsfxsize=3cm \epsfysize=3cm \epsffile{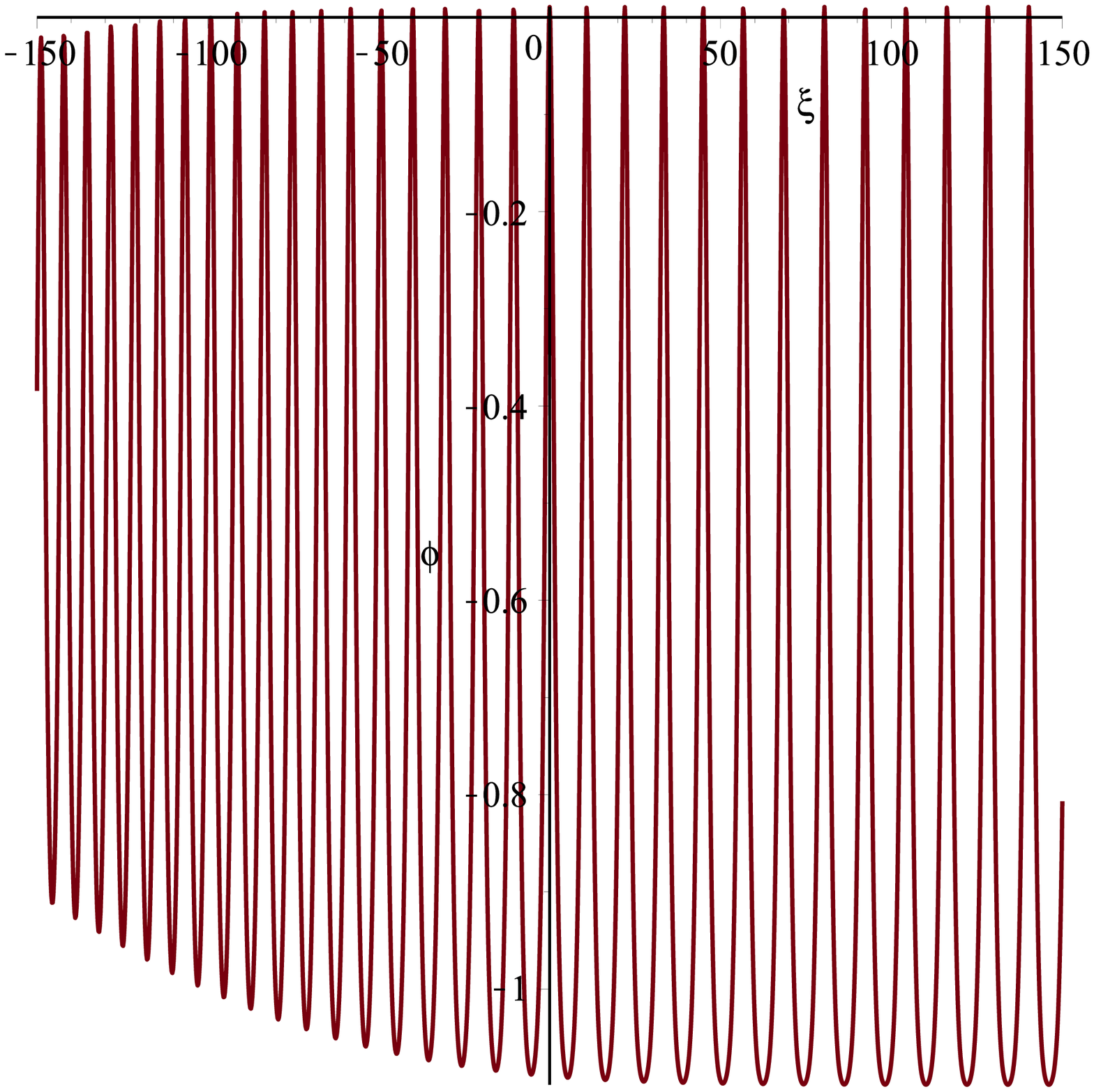}\\
\footnotesize{(a1)\ $c=0.2660295689+0.0001$ } \footnotesize{(a2)\ $(\xi,\phi)$ }\\
\footnotesize {(a)\ Initial value $(\phi(0),y(0))=(\phi_r-10^{-4},0)$}
\end{tabular}
&
\begin{tabular}{cc}
\epsfxsize=3cm \epsfysize=3cm \epsffile{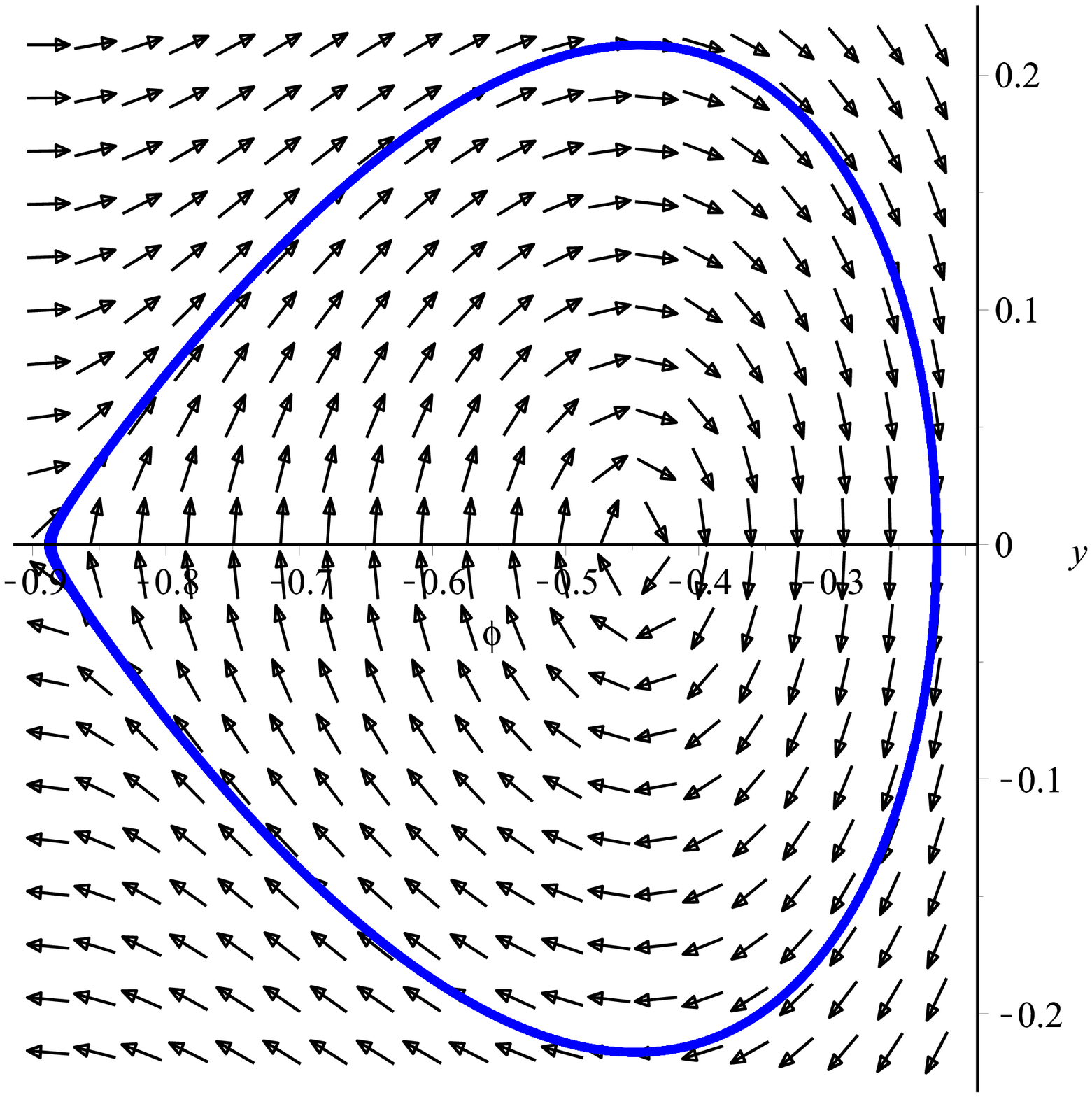}\ \ \  \epsfxsize=3cm \epsfysize=3cm \epsffile{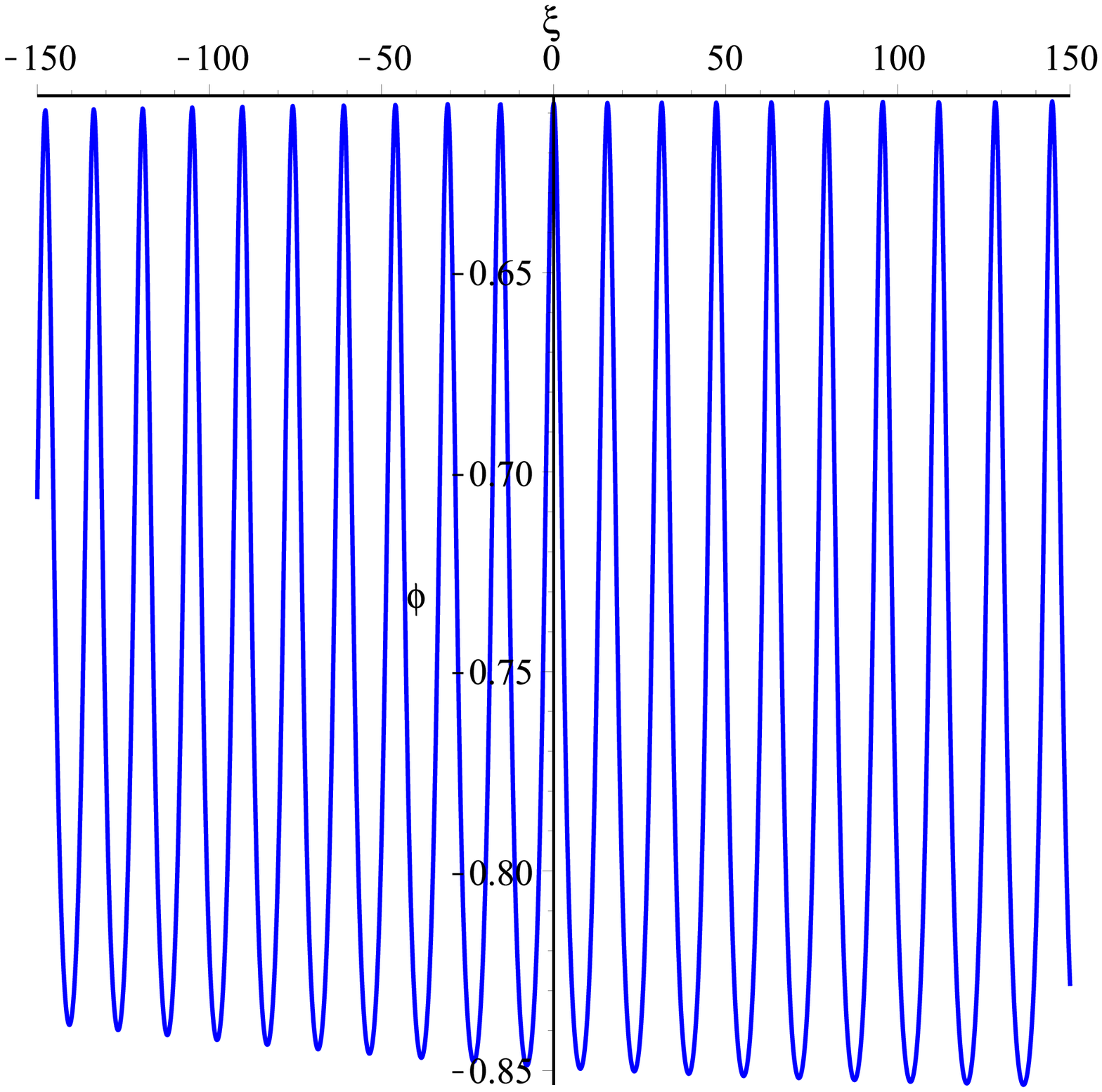}\\
\footnotesize{(b1)\ $c=0.3286393802+0.0001$ }  \footnotesize{(b2)\ $(\xi,\phi)$ }\\
\footnotesize {(b)\ Initial value $(\phi(0),y(0))=(\phi_r-10^{-4},0)$}
\end{tabular}\\
\hline  $g=0$ &
\begin{tabular}{cc}
\epsfxsize=3cm \epsfysize=3cm \epsffile{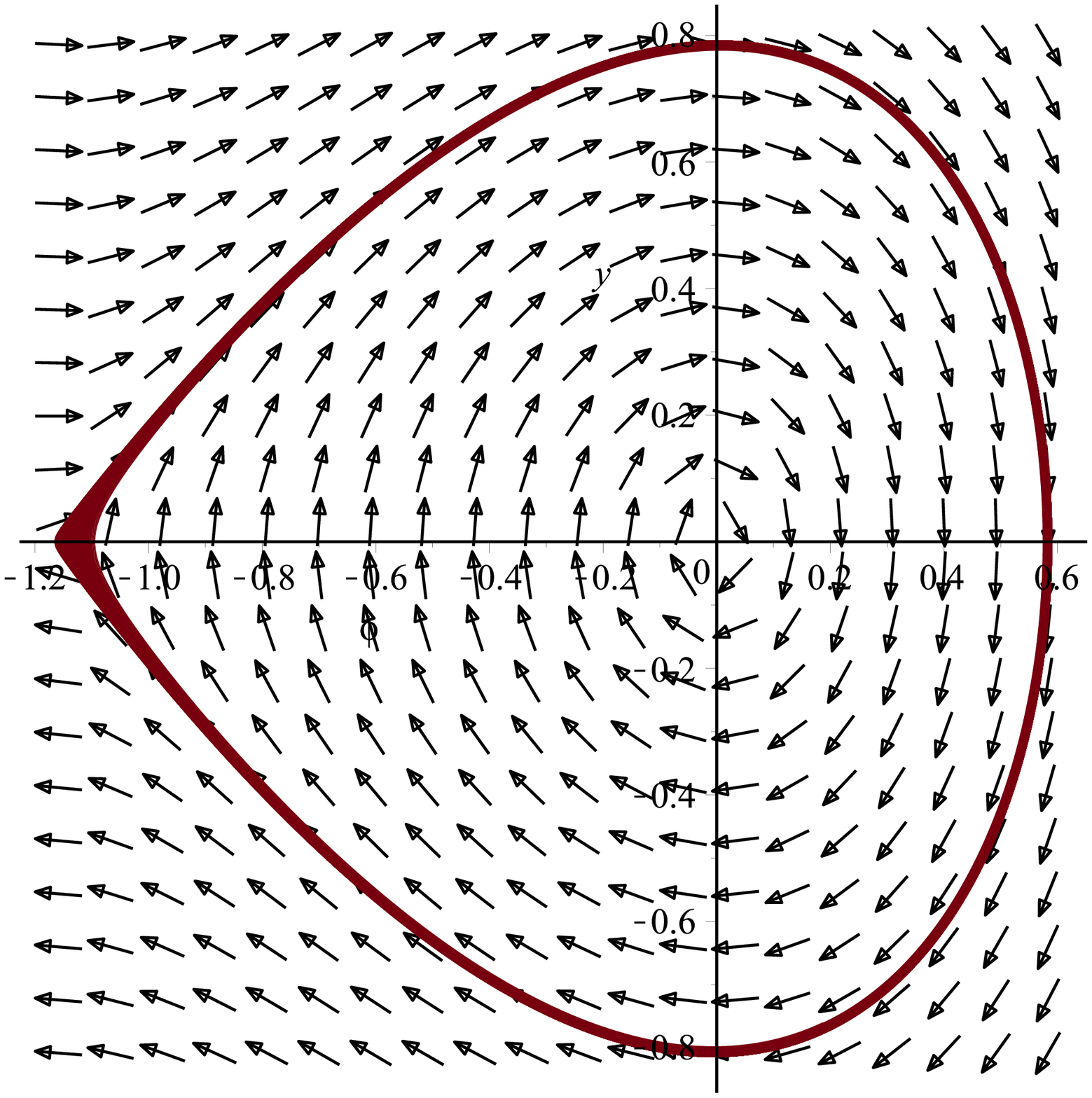}\ \ \ \epsfxsize=3cm \epsfysize=3cm \epsffile{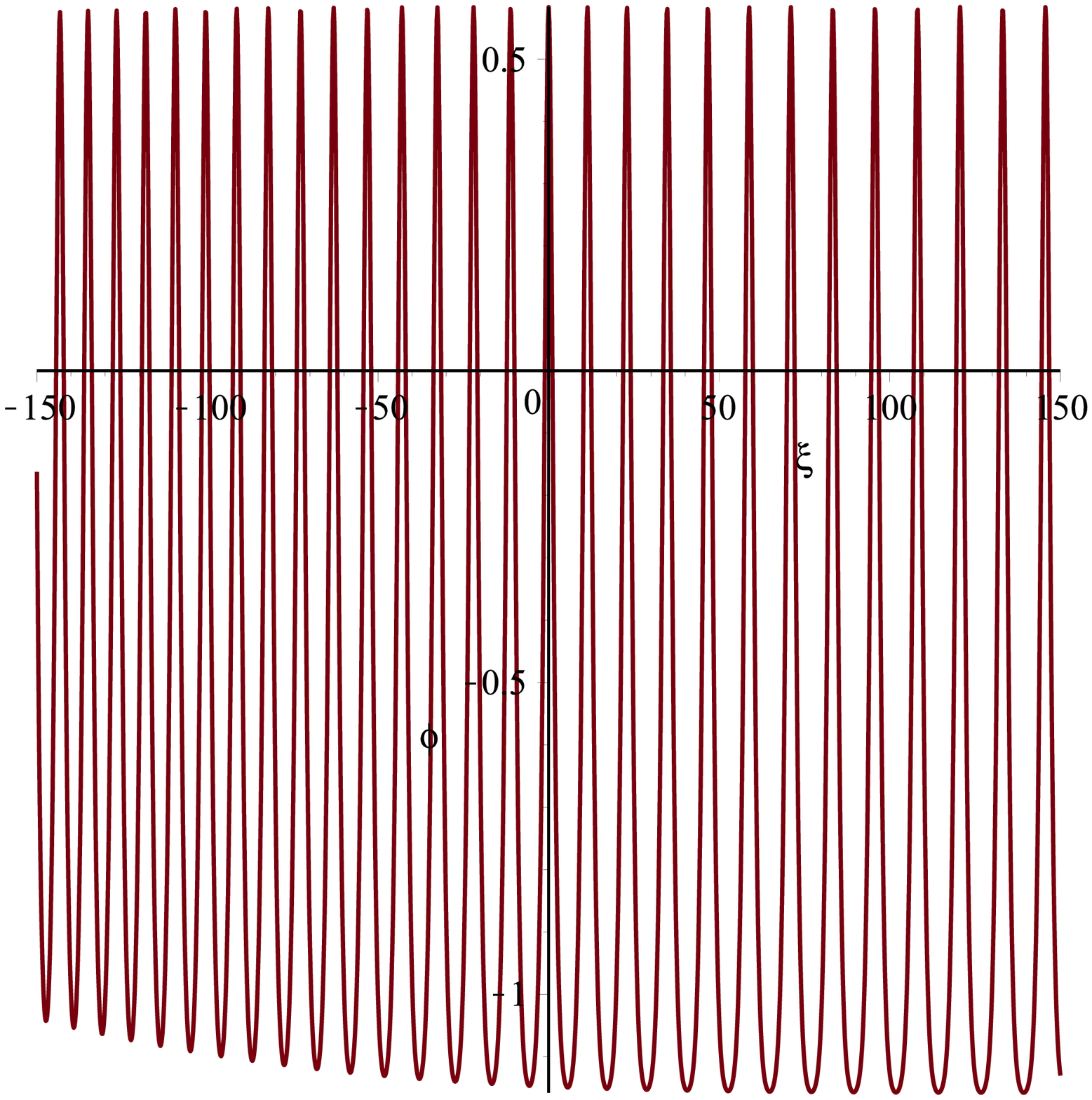}\\
\footnotesize{(c1)\ $c=0.4166666667+0.0001$ }  \footnotesize{ (c2)\ $(\xi,\phi)$ }\\
\footnotesize {(c)\ Initial value $(\phi(0),y(0))=(\phi_r-10^{-4},0)$}
\end{tabular}
&
\begin{tabular}{cc}
\epsfxsize=3cm \epsfysize=3cm \epsffile{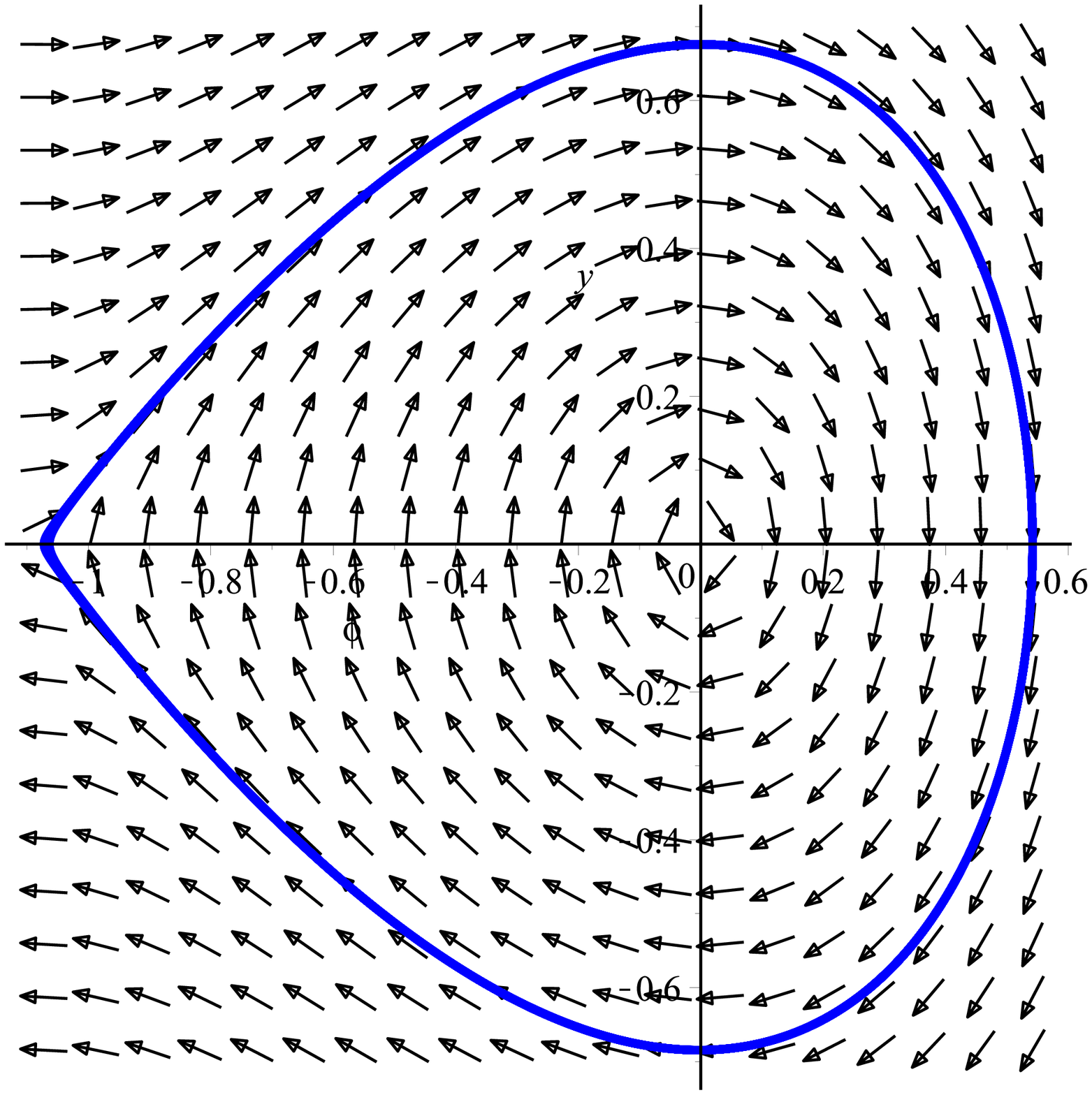}\ \ \  \epsfxsize=3cm \epsfysize=3cm \epsffile{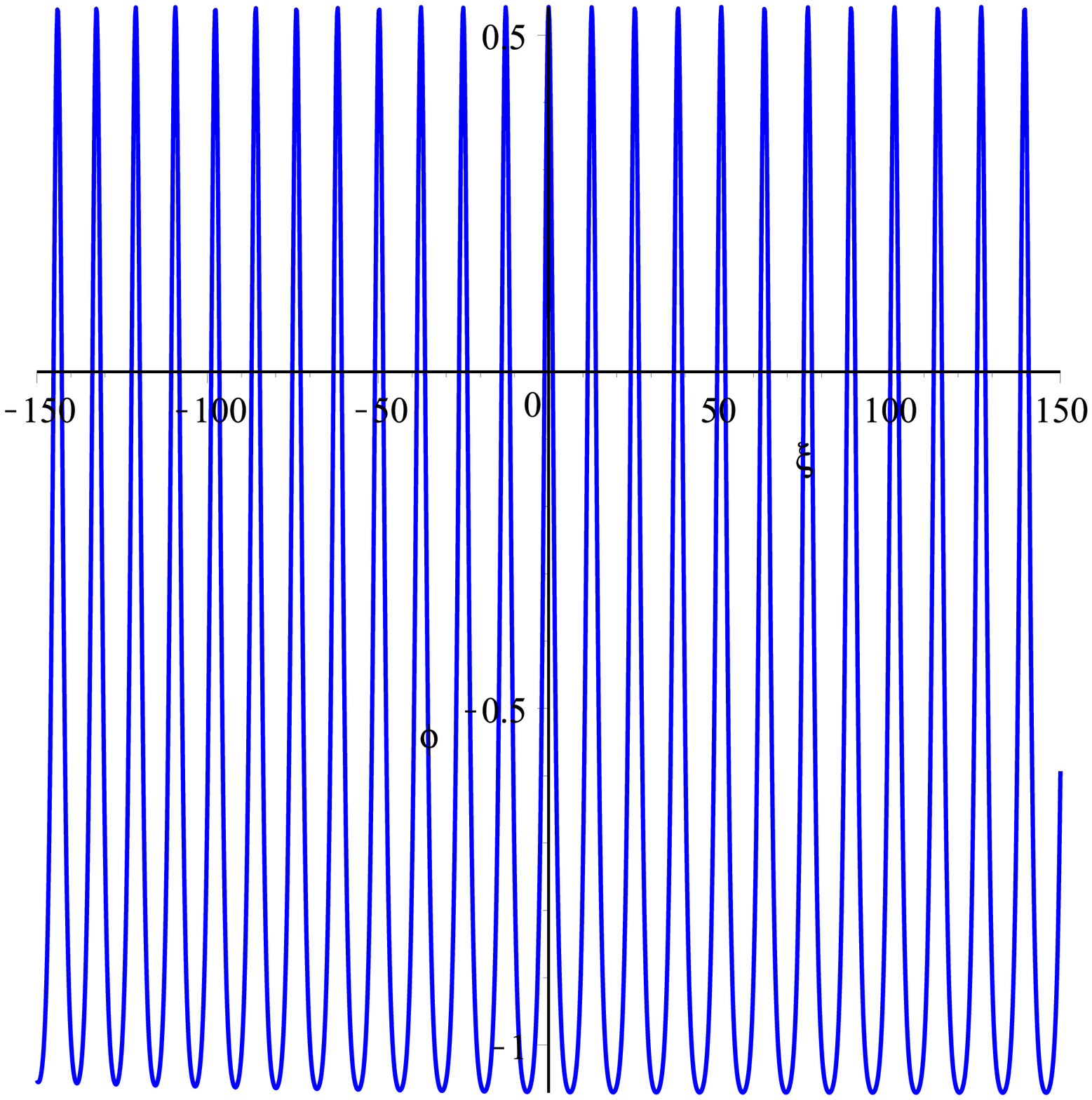}\\
\footnotesize{(d1)\ $c=0.4580398915+0.0001$ }  \footnotesize{(d2)\ $(\xi,\phi)$ }\\
\footnotesize {(d)\ Initial value $(\phi(0),y(0))=(\phi_r-10^{-4},0)$}
\end{tabular}\\
\hline $g=2$ &
\begin{tabular}{cc}
\epsfxsize=3cm \epsfysize=3cm \epsffile{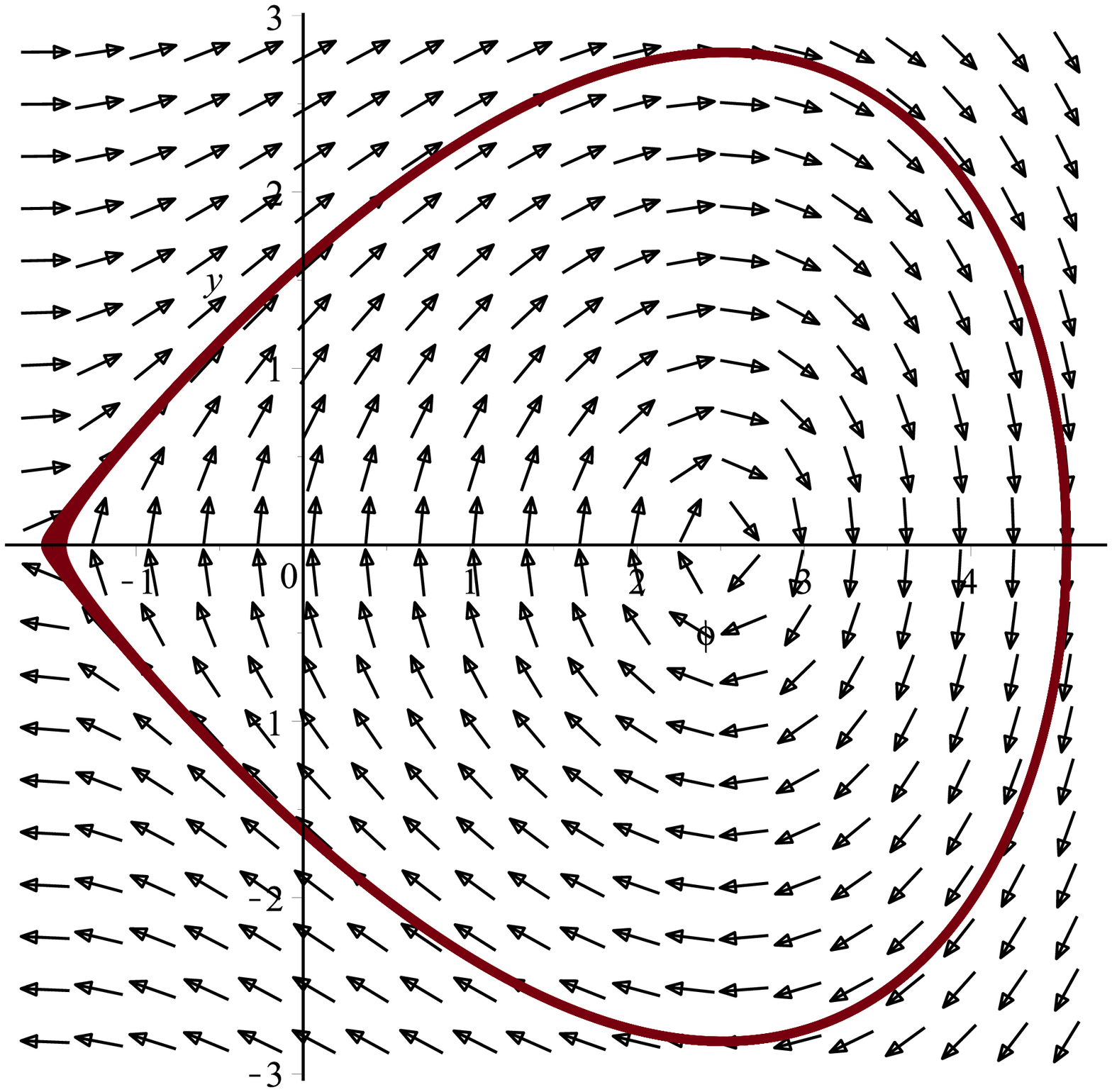}\ \ \  \epsfxsize=3cm \epsfysize=3cm \epsffile{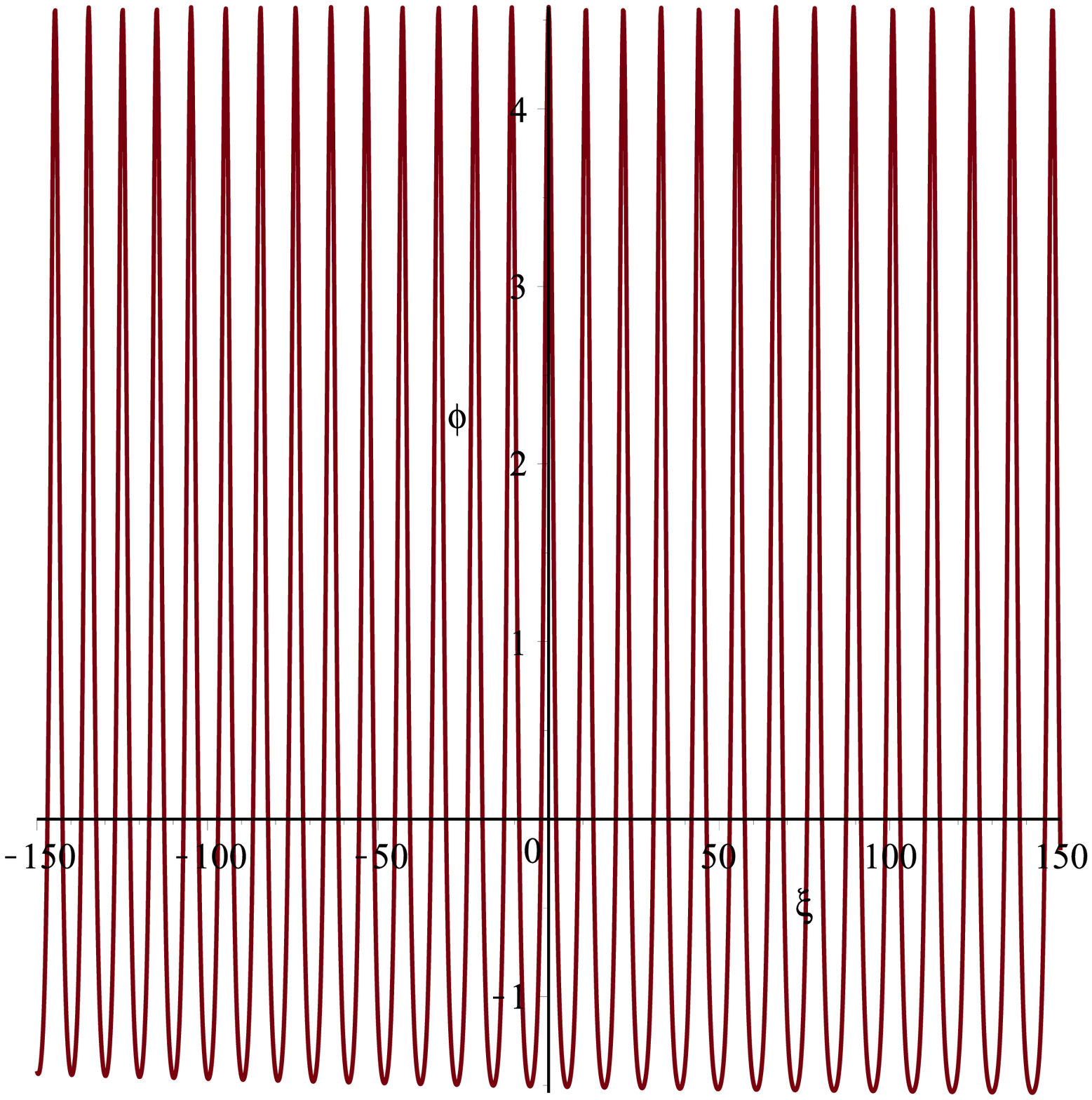}\\
\footnotesize{(e1)\ $c=1.466998871+0.0001$ }  \footnotesize{(e2)\ $(\xi,\phi)$ }\\
\footnotesize {(e)\ Initial value $(\phi(0),y(0))=(\phi_r-10^{-3},0)$}
\end{tabular}
&
\begin{tabular}{cc}
\epsfxsize=3cm \epsfysize=3cm \epsffile{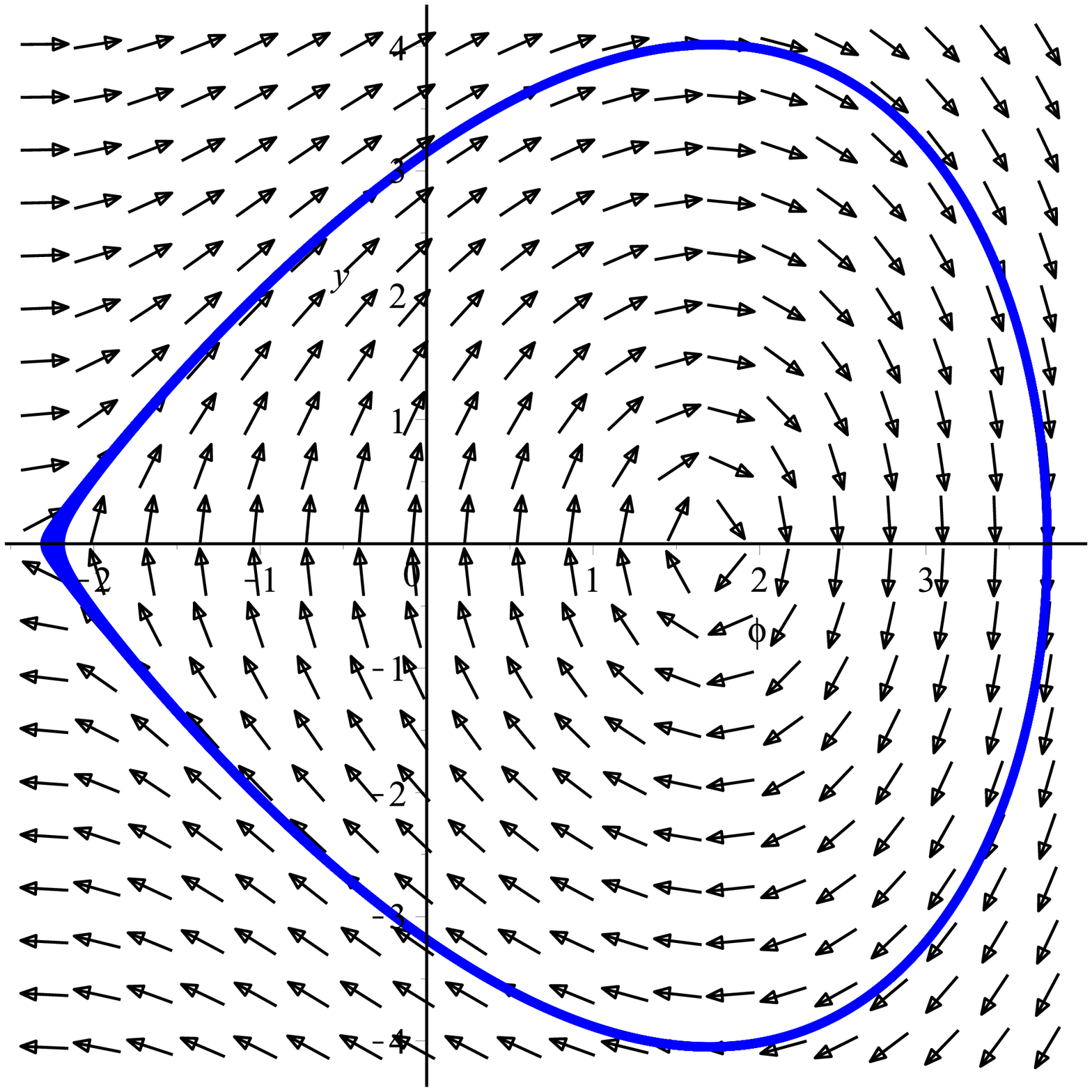}\ \ \ \epsfxsize=3cm \epsfysize=3cm \epsffile{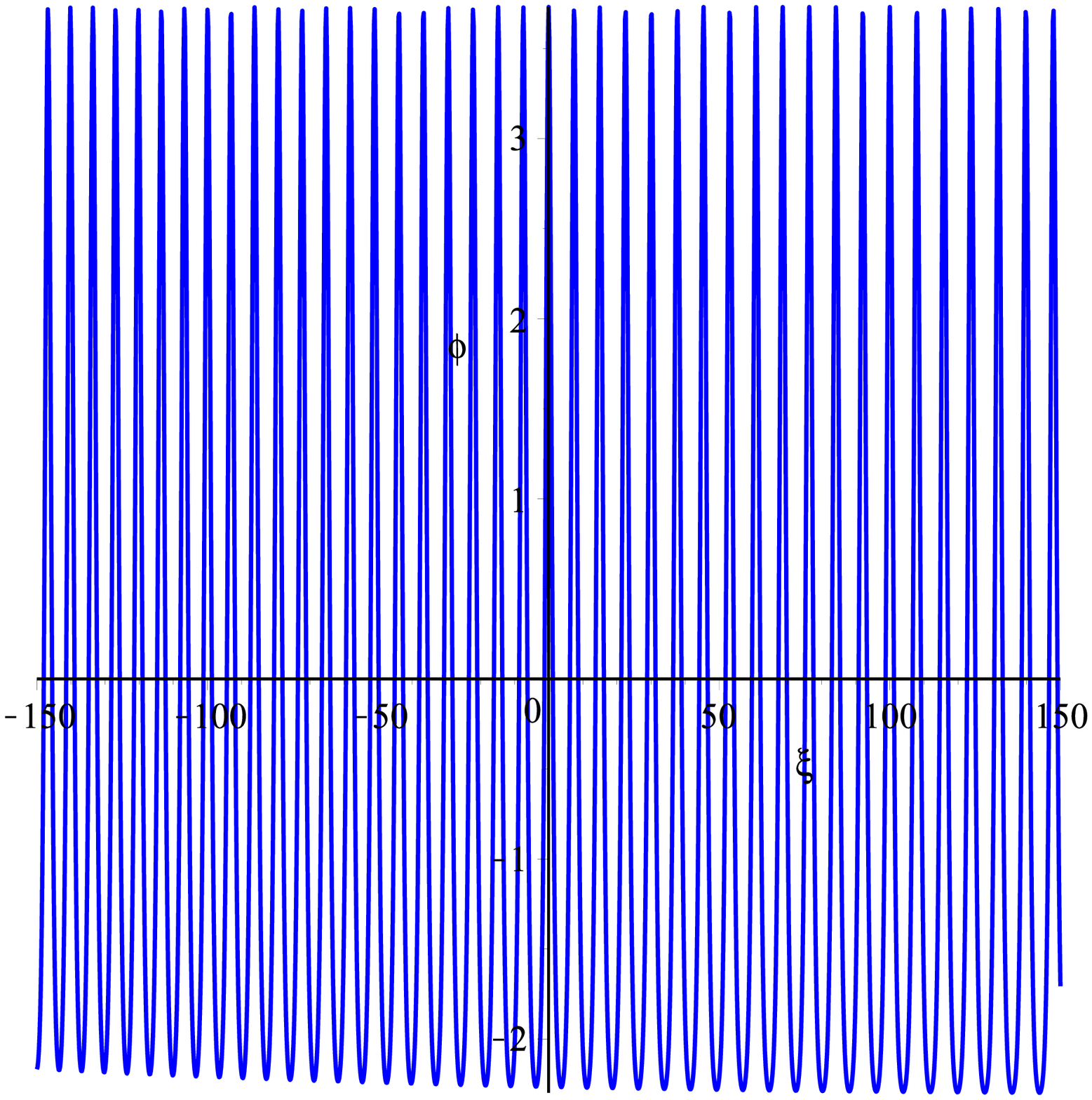}\\
\footnotesize{(f1)\ $c=0.6801960271+0.0001$ }  \footnotesize{(f2)\ $(\xi,\phi)$ }\\
\footnotesize {(f)\ Initial value $(\phi(0),y(0))=(\phi_r-10^{-3},0)$}
\end{tabular}\\
\hline
\end{tabular}
\end{center}
\end{table}
\begin{remark}
%\noindent {\bf Conclusions:}\\
The above simulations and the figures confirm the following facts: \\
1. Indeed, the solitary wave persists under small perturbations, whether the perturbations are KS or ME.\\
2. The different perturbations do affect the  proper wave speed $c$ ensuring the persistence of the solitary waves. The essential reason is that the nonlinear term $(uu_x)_x$ induces different representations of Melnikov functions. Consequently, solitary wave solutions of corresponding system exist at different proper wave speed $c$ for different kinds of small perturbations.
\end{remark}

%\noindent {\em Discussion:} As we see,  {RAJ-SIAM}

%\section{Conclusion}

\section{Acknowledgement}
{The authors would like to thank Professor A.J. Roberts for his valuable discussions when they are preparing the manuscript.} This work was jointly supported by the National Natural Science Foundation of China under Grant (No. 11671176, 11931016), Natural Science Foundation of Fujian Province under Grant (No. 2021J011148), Fujian Province Young Middle-Aged teachers education scientific research project (No. JAT210454, No. JAT200670) and Teacher and Student Scientific Team Fund of Wuyi University (Grant No. 2020-SSTD-003).

\section{ Conflict of Interest}
 The authors declare that they have no conflict of interest.

\section{Data Availability Statement}

No data was used for the research described in the article.

\section{Author contributions}

All the authors have the same contributions to the paper.

%\end{multicols}
\end{document}